%% file: ordinaryNonUniruled.tex
\documentclass[twoside,10pt]{amsart}
  \usepackage{mabliautoref}
  \usepackage{amsmath,amsfonts,amssymb,amsthm,mathtools,multirow,tabls,mathrsfs}

  \usepackage[a4paper,margin=1in]{geometry}  
  \usepackage{lmodern}
  \usepackage[all]{xy}
  \usepackage{scalerel}
  \usepackage[utf8]{inputenc}
  \usepackage[T1]{fontenc}
  \usepackage[shortlabels]{enumitem}
  \usepackage{cancel}
  
  \newcommand{\eqdef}{\overset{\mathrm{def}}{=\joinrel=}}

\let\counterwithin\relax
\usepackage{chngcntr}

\numberwithin{equation}{section}

\input{macrosMath}

\input{includeNice3}
  
  \usetikzlibrary{tikzmark,shapes.geometric, arrows, positioning, fit, calc,decorations.pathmorphing,shapes}
\tikzstyle{e} = [ellipse, minimum width=1cm, minimum height=0.5cm,text centered, draw=black]
\tikzstyle{arrow} = [thick,-,>=stealth]
\tikzstyle{connection}=[inner sep=0,outer sep=0]

  \author[Zs. Patakfalvi]{Zsolt Patakfalvi}
  \address{\'Ecole Polytechnique F\'ed\'erale de Lausanne, Chair of Algebraic Geometry \newline 
    \indent MA C3 625 (Bâtiment MA), Station 8, CH-1015 Lausanne}
  \email{zsolt.patakfalvi@epfl.ch}

  \author[M. Zdanowicz]{Maciej Zdanowicz}
  \address{\'Ecole Polytechnique F\'ed\'erale de Lausanne, Chair of Algebraic Geometry \newline 
    \indent MA C3 585 (Bâtiment MA), Station 8, CH-1015 Lausanne}
  \email{maciej.zdanowicz@epfl.ch}

  \title[Ordinary Calabi-Yau varieties are not uniruled]{Ordinary varieties with trivial canonical bundle are not uniruled}
  \date{\today}

    \DeclareMathOperator{\AR}{AR}
    
\renewcommand{\QQ}{\mathbb{Q}}
\renewcommand{\FF}{\mathbb{F}}

\begin{document}

\begin{abstract}
We prove that smooth, projective, $K$-trivial, weakly ordinary varieties over a perfect field of characteristic $p>0$ are not geometrically uniruled.  We also show a singular version of our theorem, which is sharp in multiple aspects.  Our work, together with Langer's results, implies that varieties of the above type have strongly semistable tangent bundles with respect to any polarization. 
\end{abstract}

\maketitle


\input{introduction.tex}

\input{preliminaries.tex}

\input{Witt.tex}

\input{top_Witt.tex}


\input{CY.tex}


\bibliographystyle{amsalpha} 
\bibliography{includeNice}
\end{document}

%% file: macrosMath.tex
\newcommand{\cE}{\mathscr E}
\newcommand{\cF}{\mathscr F}
\newcommand{\cG}{\mathscr G}
\newcommand{\cH}{\mathscr H}

\newcommand{\cN}{\mathscr N}
\newcommand{\cO}{\mathscr O}

\newcommand{\cT}{\mathscr T}

\newcommand{\scr}{\mathcal}

\newcommand{\sE}{\scr{E}}
\newcommand{\sF}{\scr{F}}
\newcommand{\sG}{\scr{G}}

\newcommand{\sO}{\scr{O}}

\newcommand{\sQ}{\scr{Q}}

\newcommand{\sT}{\scr{T}}

\newcommand{\bb}[1]{\mathbb{#1}} 
\newcommand{\FF}{\bb{F}}

\newcommand{\NN}{\bb{N}}
\newcommand{\PP}{\bb{P}}
\newcommand{\QQ}{\bb{Q}}

\newcommand{\ZZ}{\mathbb{Z}}

\renewcommand{\phi}{\varphi}
\newcommand{\isom}{\simeq} 
\newcommand{\piet}{\pi^{\acute{e}t}_{1}}
\newcommand{\et}{{\rm\acute{e}t}}
\newcommand{\ra}{\longrightarrow}    
\newcommand{\isomto}{\xrightarrow{\,\smash{\raisebox{-0.65ex}{\ensuremath{\scriptstyle\sim}}}\,}}

\DeclareMathOperator{\btor}{b-tors}

\DeclareMathOperator{\Ext}{Ext}

\DeclareMathOperator{\Hom}{Hom}

\DeclareMathOperator{\Spec}{Spec}

\DeclareMathOperator{\Ab}{Ab}

\newcommand{\cExt}{{\mathscr E}\kern -.5pt xt}
\newcommand{\cHom}{\mathscr{H}\kern -.5pt om}

\newcommand{\stacksproj}[1]{{\cite[Tag~\href{http://stacks.math.columbia.edu/tag/#1}{#1}]{stacks-project}}}



%% file: introduction.tex

\section{Introduction}


{\scshape We work over a perfect field $k$ of characteristic $p>0$.} 
Consider a smooth projective $K$-trivial variety $X$ of dimension $d$, where $K$-trivial means that  $K_X \sim 0$. The similitude or the dissimilarity of $X$ to a $K$-trivial variety of characteristic zero can be measured both geometrically and cohomologically.  On the geometric side, this is done typically via properties exhibiting behaviors that do not happen in characteristic zero, and which hold usually over loci that are not dense in moduli. Some examples are:
\begin{enumerate}[leftmargin=43pt]
\item[\label=(unirat)] \label{itm:unirational} $X$ is \emph{unirational}: it admits a generically finite dominant rational map $\bP^n \dashrightarrow X$,
\item[(unirul)] \label{itm:uniruled} $X$ is \emph{uniruled}: it admits a generically finite dominant rational map $\bP^1 \times Z \dashrightarrow X$,
\item[(rcc)] \label{itm:rationally_chain_connected} $X$ is \emph{rationally chain connected}: two general points of $X$ can be connected by a chain of rational curves.
\end{enumerate}
Below we also list some of the properties used on the cohomological side. The  first four of these are genericity conditions that typically hold over dense open sets of moduli spaces, while last one typically holds on complements of such open dense sets:
\begin{enumerate}[leftmargin=43pt]
\item[(w-ord)] \label{itm:weakly_ordinary} $X$ is \emph{weakly ordinary}: the action of the Frobenius morphism on $H^{d}(X, \sO_X)$ is bijective,
 \item[(ord)] \label{itm:ordinary} $X$ is \emph{ordinary}, see \cite[Definition 7.2]{Bloch_Kato_p_adic_etale_cohomology},
\item[(Witt)] \label{itm:Witt} $ H^d\left(X, W\sO_{X, \bQ} \right) \neq 0$,
\item[(Br)] \label{itm:Brauer} the formal Brauer group \big(the deformation space of $0 \in H^2( \bG_m)$\big)  is isomorphic to $\wh{\bG_m}^1$, 
\item[(Shio)] \label{itm:cycles} $H^2(X, \sO_X) \neq 0$, but algebraic cycles span the entire $H^2_{\et}\left(X_{\ok}, \bQ_l(1)\right)$ 
(Shioda supersingular).
\end{enumerate} 
%
%
%
Many important results and conjectures of the area concern connections between the above notions. In the following diagram, where the top row contains the geometric and the bottom row the cohomological notions, we portray some of what is known for
 \fbox{$X$ smooth projective over $k$ with $K_X \sim 0$ and $d:= \dim X$}:
 \begin{equation}
 \label{eq:geom_coho}
 \xymatrix@C=60pt@R=55pt{
 &  &  \textrm{not (unirul)} 
  \ar@{=>}[r]
  & \textrm{not (rcc)}  
  \ar@{=>}[r]
  & \textrm{not (unirat)}
  \ar@/^1pc/@{==>}[l]^{\textrm{if $X=$ surface}}
\\
 %
 \textrm{(Br)} & 
 \ar@{<==>}[l]^{\textrm{if $X=$K3 surface}}
 \textrm{(ord)} \ar@{=>}[r] & 
 \ar@/^1pc/@{==>}[l]^{\textrm{if $X=$K3 surface or abelian variety}}
 \textrm{(w-ord)} 
 \ar@[red]@{=>}[u]^{\textrm{ \autoref{thm:main_thm_intro_smooth} }}
 \ar@{=>}[r]^{\parbox{70pt}{\scriptsize if $H^{d -1}(X, \sO_X)=0$, by \autoref{lem:wo_o_cohomology_ordinary}}}
 \ar@/_1pc/@{=>}[r]_{\parbox{70pt}{\scriptsize \ \\[3pt] if $H^{d -1}(X, \sO_X)\neq 0$, by \autoref{example:singular_enriques_surfaces}}}|(0.55){\raisebox{5pt}{\scaleobj{1.3}{\xcancel{\ }}}}
 & 
 \textrm{(Witt)} 
  \ar@/_1.7pc/@{=>}[ul]|(0.35){\qquad \parbox{75pt}{\scriptsize \cite{Esnault_Varieties_over_a_finite_field_with_trivial_Chow_group_of_0_cycles_have_a_rational_point}, see \autoref{lem:uniruled_witt_vector} for  alternative proof} }
 \ar@{<==>}[r]_{\parbox{80pt}{\scriptsize{if $X=$ K3 surface \\ \cite[Theorem 4.1]{Artin_Mazur} }}} 
 & 
 %
 %
 \textrm{not (Shio)} 
 \ar@{==>}[u]|(0.3){\parbox{50pt}{ \scriptsize if $X=$  K3 surface \cite{Shioda_An_example_of_unirational_varieties_in_characteristic_p}} } 
 }
 \end{equation}
%
%
%
It is particularly interesting to connect the cohomological and the geometric side. In particular, this is the topic of our main result, also shown by the red arrow on  diagram \autoref{eq:geom_coho}. For historical comments and additional information on diagram \autoref{eq:geom_coho}, see \autoref{ss:historical_remarks}.

\begin{theorem}[Theorem~\ref{thm:weakly_ordinary_cy_not_uniruled}]\label{thm:main_thm_intro_smooth}
{\scshape Smooth case:} If $X$ is a smooth weakly ordinary projective variety over $k$ satisfying $K_X \sim 0$, then $X$ is not geometrically uniruled. 
\end{theorem}



\noindent
We note that \autoref{thm:main_thm_intro_smooth} also works for $K_X \equiv 0$ by replacing weak ordinarity with global $F$-splitting, see \autoref{thm:weakly_ordinary_cy_not_uniruled}.

The classical works \cite{Shioda_An_example_of_unirational_varieties_in_characteristic_p,Artin_Supersingular_K3_surfaces}, shown by the right vertical arrow of \autoref{eq:geom_coho} and reviewed in \autoref{ss:historical_remarks}, are very specific to K3 surfaces.  Hence, our approach is more similar in spirit to the  diagonal arrow of \autoref{eq:geom_coho}: 
we use Witt cohomology and $p$-adic cohomology invariants. 
However, as mentioned already in the bottom row of the diagram \autoref{eq:geom_coho}, the direct connection between weak ordinarity and the non-vanishing of $H^d \left(X, W\sO_{X, \bQ} \right) $ works only if $H^{d-1}(X, \sO_X) \neq 0$ (though, in \cite[Theorem 3.2.1]{Joshi_Rajan} Joshi and Rajan proved finiteness of $H^d (X, W\sO_{X})$ over $W(k)$ in our situation). As for example $X$ could have non-trivial Albanese morphism, in which case $H^{d-1}(X, \sO_X) \cong H^1(X, \sO_X) \neq 0$ would hold, there is no hope to apply this connection directly to $X$. We solve this issue by passing to the geometric generic fiber of the maximal rationally chain connected fibration (MRCC fibration for short). In particular, most of the article is about finding a singularity class which satisfies both of the following:
\begin{itemize}
\item the geometric generic fiber of the MRCC fibration has these types of singularities, and
\item we are able to show the non uniruledness in the presence of  this class of singularities.
\end{itemize}
In other words, it is inherent for our methods (see \autoref{ss:strategy} for details) that we prove a singular version of \autoref{thm:main_thm_intro_smooth} too: 

\begin{theorem}[\autoref{thm:weakly_ordinary_cy_not_uniruled}]\label{thm:main_thm_intro_singular}
{\scshape Singular case:}
Let $X$ be a normal, $S_3$, projective, globally $F$-split variety  over $k$ with  W$\sO$-rational singularities and with $K_X\sim 0$.   Then $X$ is not geometrically uniruled. 
\end{theorem}

In \autoref{thm:weakly_ordinary_cy_not_uniruled}, we use the condition of $F$-splitting which is equivalent to weak ordinarity if $K_X$ is actually trivial (see \autoref{prop:ordinary_is_fsplit}).  This property is better suited for the purpose of our proof -- it naturally descends to the general fibers of fibrations.  Unfortunately, we were not able to prove that same statement for  W$\sO$-rationality, which explains the presence of $p$-adic cohomology in the argument.  More precisely, using proper base change, to the best of our knowledge unavailable for W$\sO$-cohomology, we managed to prove  that the general fiber of a fibration of a W$\sO$-rational variety admits a quasi-resolution  preserving $p$-adic cohomology (\autoref{prop:geom_generic_fibres_QQp_rational}, see \autoref{defin:quasi_resolution} for the definition of a quasi-resolution).  This leaves a natural question:

\begin{question}
Let $X$ be a W$\sO$-rational variety.  Assume that $f \colon X \to S$ is a fibration with geometric generic fibre reduced.  Is the geometric generic fiber also W$\sO$-rational?
\end{question}

\noindent We note that a characteristic zero analogue of the above statement, involving the notion of rationality, can in fact be proven using a simple base change argument.

\begin{remark}
The assumption of W$\sO$-rationality in \autoref{thm:main_thm_intro_singular} is satisfied in a few important special cases.  For example, if $X$ is rational, that is, there exists a resolution of singularities $f \colon Y \to X$ such that $Rf_*\sO_Y = \sO_X$, or if $X$ is a klt threefold (see \cite[Theorem 1.4]{Gongyo_Nakamura_Tanaka_Rational_points_on_log_Fano_threefolds_over_a_finite_field}).  
\end{remark}






We finish this part of the introduction by noting that by the previous work of the authors \cite{Patakfalvi_Zdanowicz_On_the_Beauville_bogomolov_decomposition}, a Beauville--Bogomolov type decomposition holds for $K$-trivial weakly ordinary varieties.  So, in the spirit of characteristic zero approach for the Beauville--Bogomolov decomposition, using recent result of Langer \cite[Corollary 3.3]{Langer_Generic_positivity_and_foliations_in_positive_characteristic} we obtain the following:
\begin{corollary}
If $X$ is a smooth projective weakly ordinary (resp. $F$-split) variety with $K_X \sim 0$ (resp. $K_X \equiv 0$), then the tangent sheaf $ \sT_X$ is strongly $H$-semistable for any ample divisor $H$ on $X$. 
\end{corollary}


\subsection{Strategy of the proof}
\label{ss:strategy}

In this section, we present a brief description of the proof of our main theorem.  Let $X$ be as in the statement of \autoref{thm:main_thm_intro_smooth}, and set $d := \dim X$.  In order to explain the argument let us first list two basic observations:
\begin{enumerate}
\item \label{itm:uniruled_trivial_top_Witt} 
if $X$ is uniruled, then both $H^d(X, W \sO_{X,\bQ})$ and $H^d_{\et}(X,\QQ_p)$ vanish by \autoref{lem:uniruled_witt_vector}, 
\item \label{itm:top_coherent_cohom_k} 
by Serre duality, as $K_X \sim 0$, we have $H^d(X, \sO_X) \cong H^0(X, \sO_X) =k$. 
\end{enumerate}
These facts motivate the initial idea to deduce that point \autoref{itm:top_coherent_cohom_k}, together with the weak ordinarity of $X$, implies that $H^d(X,W\sO_{X,\QQ}) \neq 0$, yielding a contradiction with point \autoref{itm:uniruled_trivial_top_Witt}.  The natural argument implementing this idea is based on the following series of short exact sequences of finite length Witt vector sheaves, where $V$ is the Verschiebung homomorphism:
\[
\xymatrix{
0 \ar[r] & \sO_X \ar[r]^(0.4){V} & W_{n+1}\sO_X \ar[r] & W_n\sO_X \ar[r] & 0,
}
\]
The main point is that this argument does not always work.  There are natural obstructions  given by the edge maps in the corresponding long exact sequences of cohomology.  These obstructions are called Bockstein operations, and were analyzed thoroughly by Mumford \cite[Chapter 27]{Mumford_Lectures_on_curves_on_an_algebraic_surface} in relations to the phenomena of non-reducedness of the Picard scheme in characteristic $p$ geometry.  It turns out that already in dimension two the relevant obstructions could actually be non-zero.  For instance, every non-classical ordinary Enriques surface $S$ in characteristic two is weakly ordinary and it satisfies the conditions: $K_S \sim 0$ and $H^2(S,\sO_S) \neq 0$, but $H^2(S,W\sO_{S,\QQ}) = 0$, see \autoref{example:singular_enriques_surfaces}.

We note that, weakly ordinarity is also essential. In fact, even under the assumption that all Bockstein operations vanish, the top Witt vector cohomology $H^d(X,W\sO_{X,\QQ})$ could vanish.  However, if $H^{d-1}(X, \sO_X)=0$, through the vanishing of Bockstein operations, weak ordinarity implies that the integral Witt vector cohomology $H^d(X,W\sO_X)$ is in fact a non-zero $p$-torsion free $W(k)$-module which yields the necessary non-vanishing after localization by $\bQ$, see \autoref{lem:wo_o_cohomology_ordinary}.  Since $X$ is ordinary, we also show that the $F$-isocrystal $H^d(X,W\sO_{X,\QQ})$ is generated by its Frobenius fixed points, which implies that $H^d_{\et}(X,\QQ_p)$ does not vanish either (also \autoref{lem:wo_o_cohomology_ordinary}). .


Summing up, if we manage to put ourselves in the situation where $H^{d-1}(X, \sO_X) = 0$, while keeping enough assumption so that the above argument still holds, then we get a contradiction. Indeed, in this case we have that $H^d(X, W \sO_{X,\bQ})$ and $H^d_{\et}(X,\QQ_p)$ vanish by first paragraph, and at the same time they are non-zero by the argument after the first paragraph. Below, we explain how we manage to do this, but only for $H^d_{\et}(X,\QQ_p)$.

We claim that we may replace $X$ by the general fiber of its MRCC fibration possibly introducing mild singularities but gaining rational chain connectedness (see \autoref{example:unpleasant_fibres} for an example of a variety that we need to rule out as a fiber).  First, in our situation the assumption of weak ordinarity is equivalent to $F$-splitting which is inherited by general fibers of morphisms \cite[Cor 2.5]{Gongyo_Li_Patakfalvi_Schwede_Tanaka_Zong_On_rational_connectedness_of_globally_F_regular_threefolds}.  Then in this situation normality is also inherited by \cite[Thm 1.3.(3)]{Ejiri_When_is_the_Albanese_morphism_an_algebraic_fiber_space_in_positive___characteristic?}.
As already mentioned above, in \autoref{prop:geom_generic_fibres_QQp_rational} we also show that $X$, while not necessarily W$\sO$-rational, still admits a quasi-resolution $Y \to X$ which do not change $p$-adic cohomology.  Now, by a fundamental group argument using weak ordinarity and the Artin--Schreier sequence (see \autoref{lem:fundamental_group_artin_schreier} and \autoref{thm:fundamental_group_rcc}) one sees that $H^{d-1}(X,\sO_X) = H^1(X,\sO_X) = 0$, after possibly replacing $X$ with its \'etale cover, and hence $H^d\left(X,W\sO_{X,\QQ}\right)$ and $H^d_{\et}(X,\QQ_p)$ are non-zero.  The total space $Y$ of the postulated quasi-resolution is at least uniruled and hence satisfies the vanishing $H^d_{\et}(Y,\QQ_p)$ by \autoref{lem:uniruled_witt_vector}.  This gives a contradiction and hence finishes the proof.




\subsection{Historical remarks about the K3 surface case}
\label{ss:historical_remarks}
The aforementioned peculiar property of existence of global parametrizations, that is uniruledness or unirationality, of varieties in characteristic $p>0$ of non-negative Kodaira dimension already occurs for simple examples such as hypersurfaces.  For instance, in \cite{Shioda_An_example_of_unirational_varieties_in_characteristic_p} Shioda provides a remarkable computational proof of the fact that Fermat hypersurfaces
\[
X_n = \{x^n + y^n + z^n + u^n = 0\} \subset \PP^3_k, \quad k \text{ is an algebraically closed field of characteristic } p>0
\]
are \emph{unirational} if there exists a positive integer $e$ such that $p^e \equiv -1\pmod{n}$.   In the same paper, Shioda deduces that under the same condition on $p$ and $n$, the second cohomology group of surfaces $X_n$, in the $\ell$-adic or crystalline sense, is generated by algebraic cycles. Consequently, the Picard rank $\rho_{X_n}$ is maximal possible, equal to the second Betti number.  In general surfaces satisfying this condition are called \emph{Shioda supersingular}.  

In the special case when $n = 4$ and for $p \equiv 3\pmod{4}$, by adjunction formula, the quartic surface $X_4$ is a Shioda supersingular K3 surface. In the paper \cite{Artin_Supersingular_K3_surfaces}, Artin investigated the necessary condition for a K3 surface to be Shioda supersingular using Brauer groups.  In particular, he proved that for a K3 surface $X$ the Picard rank satisfies the inequality $\rho_X \leq 22 - 2h_X$, if the height $h_X$ of the formal Brauer group $\Phi_X$ is finite.  Consequently, no K3 surfaces with finite $h_X$ could be Shioda supersingular, and hence unirational.  The number $h_X$ turns out to be equal to the rank over $W(k)$ of the module $H^2(X,W\sO_X)$, where $W\sO_X$ is the sheaf of Witt vectors of $\sO_X$, the Dieudonn\'e module associated to $\Phi_X$ \cite[Proposition 2.13]{Artin_Mazur}.  By a simple computation (see \autoref{lem:wo_o_cohomology_ordinary}) one observes that $h_X$ is equal to one in the particular case when the Frobenius action on the one-dimensional top cohomology group $H^2(X,\sO_X)$ is bijective, that is, the surface $X$ is weakly ordinary.  On the other hand, as a consequence of Tate conjectures for K3 surfaces $h_X$ being infinite is in fact equivalent to $X$ being Shioda supersingular \cite{Charles_The_Tate_conjecture_for_K3_surfaces_over_finite_fields}.
Shioda also conjectures that in this case $X$ is unirational, which has been an active area of research in the past few decades (see, e.g., \cite{Shioda_Some_results_on_unirationality_of_algebraic_surfaces,Rudakov_Shafarevich_Supersingular_K3_surfaces_over_fields_of_characteristic_two,Liedtke_Supersingular_K3_are_unirational,Lieblich_On_the_unirationality_of_supersingular_K3_surfaces,Bragg_Lieblich_Perfect_points_on_genus_one_curves_and_consequences_for_supersingular_K3_surfaces}).


\subsection{Acknowledgements}

The authors would like to thank Piotr Achinger, Hélène Esnault, Adrian Langer, Kay Rülling and Jakub Witaszek for useful conversations and remarks.  

The work on the article was partially supported by the following grants: grant \#200021/169639 from the Swiss National Science Foundation,  ERC Starting grant \#804334, and by the National Science Foundation  Grant No. DMS-1440140 while the first author was in residence at the Mathematical Sciences Research Institute in Berkeley, California, during the Spring semester of 2019. 


%% file: preliminaries.tex

\section{General preliminaries}
\label{sec:general_preliminaries}

In this preliminary section we gather some definitions and results required in the following considerations. 


\subsection{Notation and conventions}

In the present paper, unless stated otherwise, all schemes are defined over a perfect field $k$ of characteristic $p>0$.  By $W_n(k)$ (resp. $W(k)$) we denote the ring of Witt vectors of length $n$ (resp. infinite length) and by $K$ the fraction field of the ring $W(k)$.  A \emph{variety} is a separated scheme of finite type defined over $k$.  We say that a morphism $f \colon X \to Y$ is a \emph{fibration}, if it is proper and satisfies $f_*\sO_X = \sO_S$. 


\subsection{$F$-splitting and ordinarity}
\label{ss:F_splitting_ordinarity}

Let $f \colon X \to S$ be a morphism of schemes over a perfect field $k$.  We recall that the $e$-th relative Frobenius $F^e_{X/S}$ morphism is defined by the following commutative diagram:
\[
\xymatrix@C=50pt@R=13pt{
  X \ar[rd]|{F^e_{X/S}}\ar@/_1.2pc/[rdd]_{\pi} \ar@/^1.2pc/[rrd]^{F^e_X}& & \\
   & X^e_S \ar[r]\ar[d]_{\pi'}\ar@{}[dr]|-{\square} & X \ar[d]^{\pi} \\
   &  S \ar[r]_{F^e_S} & S.
}
\]

\begin{definition}
 We say that $f$ is \emph{globally $F$-split} if there exists a splitting of the natural map $\sO_{X^e_S} \to F^e_{X/S,*}\sO_X$, as a homomorphism of $\sO_{X^e_S}$-modules.  In particular, a scheme $X$ defined in characteristic $p$ is globally $F$-split if the natural map $X \to \Spec(\FF_p)$ is \emph{globally $F$-split} (equivalently, the natural homomorphism $\sO_X \to F_*\sO_X$ is globally split).
\end{definition}

\noindent Suppose $f \colon X \to S$ is a globally $F$-split morphism over $k$, and let $T \to S$ be a morphism.  Taking a base change of a splitting we see that $f_T \colon X_T \to T$ is also globally $F$-split.  In particular, fibres of a globally $F$-split morphism are globally $F$-split.  We now recall the following result.  In the proof, we adopt a few definitions included in \cite{Ejiri_When_is_the_Albanese_morphism_an_algebraic_fiber_space_in_positive___characteristic?}.  In particular, we say that a $\QQ$-divisor $\Delta$ is $\ZZ_{(p)}$-Weil (resp. $\ZZ_{(p)}$-Cartier) if there exists number $N$ coprime to $p$ such that $N \cdot \Delta$ is Weil (resp. Cartier).  We note that such conditions can always by exhibited using a number of the form $p^f - 1$, for $f \in \NN$.

\begin{proposition}
\label{prop:f_split_fibration_normal}
Let $X$ be a globally $F$-split normal projective variety over $k$ such that $K_X$ is $\bQ$-Cartier with index prime-to-$p$ and let $f \colon X \to S$ be a fibration.  Then there exists a non-empty open subset $U \subset S$ such that $f_{|U} \colon X_U \to U$ is a globally $F$-split morphism.  Moreover, there exists an open subset $V \subset X$ such that $V \to U$ is smooth and $\codim_X(X \setminus V) = 2$.  
\end{proposition}

\begin{proof}
First, using the technique of \cite[Proposition 3.12]{Schwede_Smith_Globally_F_regular_an_log_Fano_varieties} we observe that there exists a divisor $\ZZ_{(p)}$-Cartier divisor $\Delta$ such that $(X,\Delta)$ is an $F$-split pair and $K_X + \Delta \sim_{\ZZ_{(p)}} 0$.  Localizing the splitting of $X$ at the generic point $\eta$ of $S$ and then using the argument of \cite[Proposition 2.6]{Gongyo_Li_Patakfalvi_Schwede_Tanaka_Zong_On_rational_connectedness_of_globally_F_regular_threefolds} for $D = 0$, we see that the geometric generic fibre $(X_{\overline{\eta}},\Delta_{\overline{\eta}})$ is $F$-split.  By \cite[Theorem B]{Patakfalvi_Schwede_Zhang_F_singularities_in_families} the same condition holds for all geometric fibres in the neighbourhood of $\eta$.  Using \cite[Proposition 5.8 (3)]{Ejiri_When_is_the_Albanese_morphism_an_algebraic_fiber_space_in_positive___characteristic?}, we may now assume that $X \to S$ is actually locally $F$-split (as defined in \cite{Ejiri_When_is_the_Albanese_morphism_an_algebraic_fiber_space_in_positive___characteristic?}) and therefore application of \cite[Proposition 5.7]{Ejiri_When_is_the_Albanese_morphism_an_algebraic_fiber_space_in_positive___characteristic?} yields the desired result. 
\end{proof}

\begin{corollary}\label{cor:f_split_normal_fibres}
Let $X$ be a globally $F$-split normal projective variety over $k$ such that $K_X$ is $\bQ$-Cartier with index prime-to-$p$, and let $f \colon X \to S$ be a fibration.  Then there exists a non-empty open subset $U \subset S$ such that for every geometric point $u \in U$ the fibre $X_u$ is globally $F$-split and normal. 
\end{corollary}

In the following statement, we need notions of a dualizing complex $\omega^\bullet_X$ and a dualizing module $\omega_X$ of a variety $X$ defined over $k$.  For a detailed description, we refer to \stacksproj{08XG} and \stacksproj{0A85}.  We note that by definition $\omega_X \isom \cH^{-\dim X} \left( \omega_X^{\bullet} \right)$, and that for normal varieties $\omega_X$ is in fact an $S_2$ divisorial sheaf $\cO_X(K_X)$ corresponding to the canonical divisor $K_X$.

    \begin{proposition}\label{prop:duality_for_S_k}
    Let $X$ be a projective variety over $k$.  Assume that $X$ satisfies Serre's condition $S_k$.  Then  there is a functorial isomorphism $H^{d-i}(X,\omega_X) \isom H^i(X,\cO_X)^\vee$ for $i<k-1$.  Moreover, there exists a natural surjective morphism $H^{d-k+1}(X,\omega_X) \to H^{k-1}(X,\cO_X)^\vee$.
    \end{proposition}
    \begin{proof}
        First, we observe that duality theory yields a functorial isomorphism $\Ext^{-i}(\cF,\omega^\bullet_X) \isom H^i(X,\cF)^\vee$, for every coherent sheaf $\cF$ on $X$.  In the particular, for $\cF = \cO_X$ one obtains an isomorphism $H^{-i}(\omega^\bullet_X) \isom H^i(X,\cO_X)^\vee$.
        We set $d = \dim X$ and consider the distinguished triangle:
        \[
            \omega_X[d] \ra \omega^\bullet_X \ra \cE \xrightarrow{+1}.
        \]
        Taking cohomology, we obtain an exact sequence
        \[
 H^{-i-1}(X,\cE) \ra H^{-i}(X,\omega_X[d]) \isom H^{d-i}(X,\omega_X) \ra H^{-i}(X,\omega^\bullet_X) \isom H^i(X,\cO_X)^\vee \ra H^{-i}(X,\cE)
        \]
        However, noting that $\cH^{-q}(\cE) = \cH^{-q}(\omega^\bullet_X)$ for by $q<d$ and using \stacksproj{0ECM} along with the spectral sequence:
        \[
          H^p(X,\cH^q(\cE)) \Rightarrow H^{p+q}(X,\cE),
        \]
        we see that $H^{-j}(X,\cE) = 0$ for $0 \leq j < k$. This finishes the proof.
    \end{proof}
    We also recall that Serre properties behave well with respect to the procedure of taking geometric generic fibres.
    \begin{proposition}\label{prop:property_S_k_inherited_by_fibres} 
      Let $f \colon X \to Y$ be a flat morphism of varieties over $k$.  Assume that $X$ satisfies Serre's condition $S_k$.  Then the geometric generic fibre of $f$ satisfies Serre's condition $S_k$.
    \end{proposition}
    \begin{proof}
      The generic fibre is $S_k$ by the definition of the property $S_k$.  Taking an extension to the geometric generic fibre is valid by \cite[Corollaire 6.4.2]{Grothendieck_Elements_de_geometrie_algebrique_IV_II} (flat extension).
    \end{proof}

Concerning the relation between global $F$-splitting and ordinarity, we need the following result.  We provide the proof for the reader's convenience and because the classical reference does not treat the singular case.  A slightly different argument also already appeared in \cite[Section 2.4]{Joshi_Rajan}.

\begin{proposition}[{\cite[Proposition 9]{Mehta_Ramanathan_Frobenius_splitting_and_cohomology_vanishing_for_Schubert_varieties}}]\label{prop:ordinary_is_fsplit}
Let $X$ be a normal projective variety over $k$ with $K_X \sim 0$.  Then $X$ is weakly ordinary if and only if it is globally $F$-split.
\end{proposition}

\begin{proof}
Let $U$ be the regular locus of $X$, and set $d:= \dim X$. Then we have:

$\sO_X \to F_* \sO_X$ splits 
$\explshift{-30pt}{\Longleftrightarrow}{$U$ is a big open set, and all sheaves are $S_2$}$
$\sO_U \to F_* \sO_U$ splits 
$\expl{\Longleftrightarrow}{duality}$
$F_* \omega_U \to \omega_U$ splits  
$\explshift{50pt}{\Longleftrightarrow}{$U$ is a big open set, and all sheaves are $S_2$}$ 
$F_* \omega_X \to \omega_X$ splits 
\\[0pt]  

\hspace{120pt} $\expl{\Longleftrightarrow}{$K_X \sim 0 \Rightarrow H^0(X, \omega_X) \cong k$}$ 
$H^0(X, \omega_X) \cong H^0(X,F_* \omega_X)  \to H^0(X, \omega_X)$ is surjective
\\[0pt]

\hspace{120pt} $\expl{\Longleftrightarrow}{$X$ is normal, hence \autoref{prop:duality_for_S_k} applies with $k=2$}$ 
$H^d(X, \sO_X) \to  H^d (X,F_*\sO_X) \expl{\cong}{$F$ is affine} H^d(X, \sO_X)$ is injective 
\\[0pt] 

\hspace{120pt} $\expl{\Longleftrightarrow}{a $p$-linear map between finite dimensional $k$-vector spaces, and $k$ is perfect}$  $H^d(X, \sO_X) \to  H^d (X,F_*\sO_X) \cong H^d(X, \sO_X)$ is bijective

\end{proof}

%% file: Witt.tex
\section{Witt vector and $p$-adic cohomology}
\label{sec:witt_cohomology_and_wo_rational_sings}
 
In this section we recall the basic properties of Witt vector and $p$-adic cohomology groups and direct images.  For the results related to Witt vector cohomology we refer to \cite{Berthelot_Bloch_Esnault_On_Witt_vector_cohomology_for_singular_varieties} and \cite{Rulling_Chatzistamatiou_Hodge_Witt_and_Witt_rational} for original accounts, and to
\cite[Section 2.5]{Gongyo_Nakamura_Tanaka_Rational_points_on_log_Fano_threefolds_over_a_finite_field} for an accessible summary.  The necessary results concerning $p$-adic cohomology can be found in the classical reference \cite[Exposé V, VI]{SGA5}.  For the sake of clarity, we recall most of the definitions and some useful properties.  


\subsection{Basics of W$\sO$-cohomology}

First, we recall that for every scheme defined over $\FF_p$ there exists a sequence of schemes $W_n(X)$, for $n \geq 1$, underlying the same topological space $X$ and with varying sheaves of rings $W_n\sO_X$ given by Witt vectors of length $n$, see \cite[II, Section 6]{Serre_Local_fields} or \cite[Section 2.5]{Gongyo_Nakamura_Tanaka_Rational_points_on_log_Fano_threefolds_over_a_finite_field},.  We define the sheaf $W\sO_X $ on $X$ as the following inverse limit taken in the category $\Ab(X)$ of sheaves of abelian groups
\[
W\sO_X \eqdef \varprojlim_n W_n\sO_X.
\] 
We note that unlike $W_n(X) = (X,W_n\sO_X)$ the pair $(X,W\sO_X)$ is not a scheme.

As will be stated in \autoref{def:wo_rational_morphism}, one of the defining properties of $W\sO$-rational singularities is that $R^i f_* W\sO_{Y,\QQ}=0$ for every integer $i>0$ and every quasi-resolution of the singularities (see \autoref{defin:quasi_resolution}). The question is however how one defines $R^i f_* W\sO_{Y,\QQ}$. The naive idea is to let $R^i f_* W\sO_{Y,\QQ}$ be defined as $(R^i f_* W\sO_{Y}) \otimes_{\bZ} \bQ$. According to the understanding of the authors, it is not known whether this definition is the correct one.  The main reason is that it is not known whether the $p^{\infty}$-torsion in $R^i f_* W \sO_{Y}$ is annihilated by $p^M$ for a single integer $M>0$ or not.  If the inverse system $\{ R^i f_* W_n \sO_X \}_{n \in \NN}$ satisfied the Mittag-Leffler condition, then there would not be $p^\infty$-torsion of infinite height by the arguments of \cite[Prop 2.10]{Berthelot_Bloch_Esnault_On_Witt_vector_cohomology_for_singular_varieties}.
However, to the best of our knowledge, it is again not known whether this Mittag-Leffler property always holds true. Hence, to avoid these problems, the present definitions define $R^i f_* W\sO_{Y,\QQ}$ so that its $p^{\infty}$-torsion is annihilated by $p^M$ for a single integer $M>0$. In particular, the definition becomes slightly cumbersome, in a sense that it uses the notion of localization of a category in a Serre subcategory.  

\begin{remark}
\label{remark:localization}
The notion of localization of a category in a full subcategory is well-explained in \stacksproj{02MN}.  For the original account we refer to \cite[Chapitre III]{Gabriel_Des_Categories_abeliennes}. 
\end{remark}

In our case, this is applied to the full subcategory $\Ab(X)_{\btor}$ of the category $\Ab(X)$ of abelian sheaves on a scheme $X$, defined by the formula:
\[
A \in {\rm Ob}\left(\Ab(X)_{\btor}\right) \iff \exists N \in \bZ \setminus 0  \ : \  N \cdot \id_A = 0.
\]
The category obtained by localizing  $\Ab(X)$ in $\Ab(X)_{\btor}$ is denoted  $\Ab(X)_{\QQ}$. The category  $\Ab(X)_{\QQ}$  is called the category of $\QQ$-localized abelian sheaves, and it is an abelian category. One intuitive description of this localization is by formal inverting of all arrows in $\Ab(X)$ whose kernel and cokernel are in $\Ab(X)_{\btor}$.  We denote the natural projection functor by $q \colon \Ab(X) \to \Ab(X)_{\QQ}$ or simply $(-)_{\QQ}$, and we note that:
\begin{itemize}
\item  $q(A) = 0$ if and only if $A \in \Ab(X)_{\btor}$ \cite[3.7.2, line 6]{Rulling_Chatzistamatiou_Hodge_Witt_and_Witt_rational},
\item $\Hom_{\Ab(X)_{\QQ}} (q(\sF) , q(\sG))= \Hom_{\Ab(X)}(q(\sF), q(\sG)) \otimes_{\bZ} \bQ$ for any $\sF, \sG \in \Ab(X)$ \cite[Prop 3.7.4]{Rulling_Chatzistamatiou_Hodge_Witt_and_Witt_rational}, and
 \item the functor $q$ is exact \stacksproj{02MN} . 
 \end{itemize}
 In particular, the above setting is used to define $W\sO_{X,\QQ}$:
 \begin{equation*}
  W\sO_{X,\QQ} \eqdef q( W\sO_X) \in\Ab(X)_{\QQ}. 
 \end{equation*}
The next goal is to define the derived pushforwards of $W\sO_{X,\QQ}$.
  By \cite[Cor 3.7.5]{Rulling_Chatzistamatiou_Hodge_Witt_and_Witt_rational}, if $f \colon X \to Y$ is a morphism of varieties over $k$, the natural functor $f_* \colon \Ab(X) \to \Ab(Y)$ descends to localized categories yielding $f_* \colon \Ab(X)_{\QQ} \to \Ab(Y)_{\QQ}$.  Moreover, by \cite[Prop 3.7.7]{Rulling_Chatzistamatiou_Hodge_Witt_and_Witt_rational}, both functors admit total right derived functors $Rf_*$ which satisfy compatibility condition expressed by commutativity of the diagram:
\[
\xymatrix{
D^+(\Ab(X)) \ar[d]_{D^+(q)} \ar[r]^{Rf_*} & D^+(\Ab(Y)) \ar[d]^{D^+(q)} \\
D^+(\Ab(X)_{\QQ}) \ar[r]^{Rf_*} & D^+(\Ab(Y)_{\QQ}).
}
\]
As a consequence, the cohomology groups satisfy the equation $(R^if_*A)_{\QQ} \isom R^if_*(A_{\QQ})$, and hence this setup allows ups to define $R f_* W \sO_{X,\bQ}$:

\begin{definition}
Let $f \colon X \to Y$ be a morphism of varieties over $k$.  We define the $i$-th derived pushforward of $W\sO_{X,\bQ}$ by $R^if_* \left( W\sO_{X,\QQ} \right)  \in \Ab(Y)_{\QQ}$.  In particular, the Witt vector cohomology $H^i(X,W\sO_{X,\QQ}) \in {\rm Mod(W(k))}_{\QQ}$ of a scheme $X$ defined over a perfect field $k$ is given by the formula:
\[
H^i(X,W\sO_{X,\QQ}) = R^i\Gamma  \left( W\sO_{X,\QQ} \right),
\]
where $\Gamma \colon X \to \Spec(k)$ is the structure map.
\end{definition}

\begin{remark}
\label{rem:Q_Witt_coho_proper_case}
According to \cite[Remark 2.17.(5)]{Gongyo_Nakamura_Tanaka_Rational_points_on_log_Fano_threefolds_over_a_finite_field},
if $X$ is proper over $k$, then $H^i(X,W\sO_{X,\QQ})=0$ if and only if $H^i(X,W\sO_{X}) \otimes_{\bZ} \bQ=0$. As explained in the first paragraph of the section, it is not known if the corresponding statement for a proper birational morphism holds true.
\end{remark}

Witt vector cohomology relates to finite length Witt vector cohomology via the following results. These results are deduced from Grothendieck spectral sequence for the composition of derived functor $R\lim \circ Rf_* = Rf_* \circ R\lim$ applied for the system $n \mapsto W_n\sO_X$ which satisfies Mittag-Leffler condition.

\begin{proposition}\label{prop:witt_vectors_and_lim}
\cite[Lem 2.18]{Gongyo_Nakamura_Tanaka_Rational_points_on_log_Fano_threefolds_over_a_finite_field}
Let $f \colon X \to Y$ be a morphism of varieties over $k$.  Then, for every $i \geq 0$, there exists a short exact sequence:
\[
0 \ra R^1  \varprojlim_{n} \left(R^{i-1}f_*W_n\sO_X \right) \ra R^if_* W\sO_X \ra \left( \varprojlim_n R^if_* W_n\sO_X \right) \ra 0.
\]
In particular, if $X$ is a proper scheme over  $k$, then for all $i$ the system $n \mapsto H^i(X,W_n\sO_X)$ satisfies the Mittag-Leffler condition, and hence 
\[
H^i(X,W\sO_{X,\QQ}) = q\left( \varprojlim_{n} H^i(X,W_n\sO_X)\right) \expl{=}{\autoref{rem:Q_Witt_coho_proper_case}} \left( \varprojlim_{n} H^i(X,W_n\sO_X)\right) \otimes_{\bZ} \bQ.
\]
\end{proposition}

We will also use the following result relating Witt vector cohomology with the rigid cohomology, that is, the substitute of crystalline cohomology in the singular setting.  Since the main application of the result is actually the comparison between Witt vector and $p$-adic cohomology (see \autoref{thm:clr_rigid_zero_padic} and \autoref{cor:p_adic_slope_zero_of_witt_vector}), and rigid cohomology can be treated as a black box for this purpose, we do not include any details, referring to the standard literature on the topic such as \cite{Berthelot_cohomologie_rigide}.  
 
\begin{theorem}[{\cite[Theorem 1.1]{Berthelot_Bloch_Esnault_On_Witt_vector_cohomology_for_singular_varieties}}]
\label{thm:bbe_rigid_witt_vector}
Let $k$ be a perfect field of characteristic $p > 0$ and let $X$ be a $k$-scheme of finite type. There exists a functorial isomorphism:
\[
H^{\bullet}_{\rm rig}(X/ K)^{<1} \isomto H^\bullet(X, W\sO_{X,\QQ}),
\]
where $K:=\Frac(W(k))$. 
\end{theorem}

\subsection{W$\sO$-singularities}

Since singular schemes naturally show up in our approach to the main theorem, and since we use a Witt-vector cohomology based invariant to show non-uniruledness (see \autoref{lem:uniruled_witt_vector}), we need a characteristic $p>0$ notion of singularities allowing for the control of Witt vector cohomology.  We will use the notion of W$\sO$-rational singularities introduced in the papers \cite{Rulling_Chatzistamatiou_Hodge_Witt_and_Witt_rational,Blickle_Esnault_Rational_singularities_and_rational_points}.

The definition uses the notion of \emph{quasi-resolution} (see \cite[Definition 4.3.1]{Rulling_Chatzistamatiou_Hodge_Witt_and_Witt_rational}).  For the sake of clarity we recall some of its details.  We say that a normal equidimensional scheme $X$ defined over a field $k$ is a \emph{finite quotient} if there exists a finite surjective morphism $Y \to X$ from a scheme $Y$ smooth over $k$.  A scheme is a \emph{topological finite quotient} if there exists a universal homeomorphism from (equiv. to) a finite quotient as defined above.  

\begin{definition}[{\cite[Definition 4.3.1]{Rulling_Chatzistamatiou_Hodge_Witt_and_Witt_rational}}]
\label{defin:quasi_resolution}
We say that a morphism between two integral varieties $f \colon X \to Y$ is a \emph{quasi-resolution} if the following conditions are satisfied:
\begin{enumerate}
    \item $X$ is a topological finite quotient,
    \item $f$ is projective, surjective, and generically finite,
    \item the extension of the function fields $k(Y) \subset k(X)$ is purely inseparable.
\end{enumerate}
\end{definition}
\noindent We note that if $Y$ is normal the last condition is actually equivalent to the Stein factorization $X' \to Y$ being a universal homeomorphism.  By the work of de Jong \cite{de_Jong_Smoothness_semi_stability_and_alterations} for every variety $X$ over a perfect field $k$ there exists a quasi-resolution with $X$ as a target (see \cite[Remark 4.3.2]{Rulling_Chatzistamatiou_Hodge_Witt_and_Witt_rational}).  More precisely, the alterations of de Jong could be chosen to be compositions of quotients of smooth varieties by finite group actions and quasi-resolutions.

\begin{lemma}
\label{lem:etale_base_change_for_quasi_resolutions}
An \'etale base change of a quasi-resolution is a quasi-resolution.
\end{lemma}

\begin{proof}
Let $X \to Y$ be a quasi-resolution and let $Y' \to Y$ be an \'etale morphism.  Take $U \to X$ to be a morphism furnishing the topological finite quotient structure on $X$.  We consider the following diagram of cartesian squares:
\[
\xymatrix{
    U' \ar[r]\ar[d]_{\text{\'etale}} & X' \ar[r]\ar[d]_{\text{\'etale}} & Y' \ar[d]_{\text{\'etale}} \\
    U \ar[r] & X \ar[r] & Y
}
\]
The morphism $U' \to U$ is \'etale and therefore $U'$ is smooth.  Consequently $U' \to X'$ yields a topological finite quotient structure on $X'$.  Using standard base change properties the map $X' \to Y'$ satisfies the conditions from \autoref{defin:quasi_resolution}, and hence $X' \to Y'$ is a quasi-resolution as desired.
\end{proof}

\begin{definition}[{\cite[Definition 4.4.1 and Proposition 4.4.6]{Rulling_Chatzistamatiou_Hodge_Witt_and_Witt_rational}}]\label{def:wo_rational_morphism}
Let $X$ be a variety over $k$.  We say that $X$ has \emph{$W\sO$-rational singularities} if and only if for every (equiv.\,some) quasi-resolution $f \colon Y \to X$ the following two conditions are satisfied:
\begin{enumerate}
\item the natural map induces an isomorphism $W\sO_{X,\QQ} \isom f_*W\sO_{Y,\QQ}$,
\item $R^if_*W\sO_{Y,\QQ} = 0$, for every $i > 0$. 
\end{enumerate}
\end{definition}

\begin{remark}\label{remark:some_remarks_on_WO_rationality}
We note that by compatibility of derived functor of $f_*$ and localization the following statements concerning the definition hold true:
\begin{enumerate}[(a)]
\item if $X$ is normal then the first condition is automatically satisfied,
\item the conditions are equivalent to the fact that the natural map in the derived category of $\QQ$-localized abelian sheaves $W\sO_{X,\QQ} \to Rf_*W\sO_{Y,\QQ}$ is an isomorphism.
\end{enumerate}
In what follows we say that a morphism satisfying the above conditions is \emph{W$\sO$-rational}.  Using this notion we may rephrase W$\sO$-rationality of a scheme $X$ by postulating the existence of a W$\sO$-rational quasi-resolution $Y \to X$.
\end{remark}

We thank Kay R\"ulling for the following result and the subsequent proof of \'etale invariance of $W\sO$-rational singularities.  

\begin{proposition}\label{prop:wo-rationality_finite_level}
Let $f \colon Y \to X$ be a proper, surjective and generically finite morphism of varieties over $k$.  Then the following conditions are equivalent:
\begin{enumerate}
\item\label{prop:wo-rationality_finite_level_it1} the morphism $f$ is W$\sO$-rational, that is, $W\sO_{X,\QQ} \to f_*W\sO_{Y,\QQ}$ is an isomorphism and $R^if_*W\sO_{Y,\QQ} = 0$, for every $i>0$,
\item\label{prop:wo-rationality_finite_level_it2} there exists an integer $M>0$ such that for every $n \geq 1$ the abelian sheaves $f_*W_n\sO_Y/W_n\sO_X$ and $R^if_*W_n\sO_X$, for $i \geq 1$, are annihilated by $p^M$.
\end{enumerate}
\end{proposition}

\begin{proof}
We first prove that \autoref{prop:wo-rationality_finite_level_it2} $\Rightarrow$ \autoref{prop:wo-rationality_finite_level_it1}.  We set $C_n = f_*W_n\sO_Y/W_n\sO_X$ and consider the following exact sequence of systems of abelian sheaves on $X$: 
\[
0 \ra \{ W_n\sO_X \}_{n \in \NN} \ra \{ f_*W_n\sO_Y \}_{n \in \NN} \ra \{ C_n \}_{n \in \NN} \ra 0.
\]
Using the long exact sequence of $R^\bullet\varprojlim$ groups, and recalling that $W_n\sO_X$ and $f_*W_n\sO_X$ satisfy Mittag-Leffler condition, we obtain:

\begin{center}
\begin{tikzpicture}[descr/.style={fill=white,inner sep=1.0pt}]
        \matrix (m) [
            matrix of math nodes,
            row sep=2em,
            column sep=2.5em,
            text height=1.5ex, text depth=0.25ex
        ]
        { 0 & \overbrace{W\sO_X}^{\varprojlim W_n\sO_X} & \overbrace{f_*W\sO_Y}^{\varprojlim f_*W_n\sO_Y} & \varprojlim C_n \\
            & \underbrace{R^1\varprojlim W_n\sO_X}_{\isom 0} & \underbrace{R^1\varprojlim f_*W_n \sO_Y}_{\isom 0} & R^1\varprojlim C_n & 0. \\ };

        \path[overlay,->, font=\scriptsize,>=latex]
        (m-1-1) edge (m-1-2)
        (m-1-2) edge (m-1-3)
        (m-1-3) edge (m-1-4)
        (m-1-4) edge[out=355,in=175] (m-2-2) 
        (m-2-2) edge (m-2-3)
        (m-2-3) edge (m-2-4)
        (m-2-4) edge (m-2-5);
\end{tikzpicture}
\end{center}

\noindent Taking images in the localized category, and hence annihilating the elements $\varprojlim C_n$ and $R^1\varprojlim C_n$, which are $p^M$-torsion, we obtain the desired statement for $i = 0$.  Then using \autoref{prop:witt_vectors_and_lim}, we have an exact sequence 
\[
0 \ra R^1\varprojlim R^{i-1}f_*W_n\sO_Y \ra R^if_*W\sO_Y \ra \varprojlim R^if_*W_n\sO_Y \ra 0.
\] 
Again all the elements in the relevant inverse systems are annihilated by $p^M$, for some fixed $M$, and therefore both the right and left elements in the above exact sequence are annihilated by $p^M$ as well.  Consequently, the middle term is $p^{2M}$-torsion and hence the implication is proven.  

We now approach \autoref{prop:wo-rationality_finite_level_it1} $\Rightarrow$ \autoref{prop:wo-rationality_finite_level_it2}.  For this purpose we apply \autoref{remark:some_remarks_on_WO_rationality} and the definition of $\Ab(X)_{\bQ}$ to see that there is an integer $N>0$ such that
%
%
$\sQ_i$ is $p^N$-torsion, where
\begin{equation}
\label{eq:wo-rationality_finite_level:def_Q}
\sQ_i:= 
\left\{
\begin{matrix}
\factor{f_* W \sO_Y}{W\sO_X} & \textrm{if }i=0 \\
R^i f_* W \sO_Y & \textrm{if }i\neq0 
\end{matrix}
\right.
\end{equation}
Consider the exact sequence
\begin{equation}
\label{eq:wo-rationality_finite_level:short_exact_sequence}
\xymatrix{
0 \ar[r] & W\sO_Y \ar[r]^{V^n} & W \sO_Y \ar[r] & W_n \sO_Y \ar[r] & 0
}
\end{equation}
The long exact sequence obtained from \autoref{eq:wo-rationality_finite_level:short_exact_sequence} by applying derived pushforwards sandwiches $R^i f_* W_n \sO_Y$ between $\sQ_i$ and $\sQ_{i+1}$. As the latter are $p^N$-torsion, we obtain that $R^i f_* W_n \sO_Y$ is $p^{2N}$ torsion for all $i \geq 0$ and all $n \geq 1$.

We have to show the same statement for $\sQ_0$ too. For that, consider the following diagram where all rows and columns are exact. To see this, start with the two leftmost non-zero columns that we already know are exact. Apply then the snake lemma first two the first three non-zero entries of these columns, and then to the second and third non-zero rows. 
\begin{equation*}
\xymatrix{
& 0 \ar[d] & 0 \ar[d] & 0 \ar[d] \\
0 \ar[r] & W\sO_X \ar[r]\ar[d]^{V^n} & f_*W\sO_Y \ar[r]\ar[d]^{V^n} & \sQ_0 \ar[d] \ar[r] & 0 \\
0 \ar[r] & W\sO_X \ar[r]\ar[d] & f_*W\sO_Y \ar[r]\ar[d] & \sQ_0 \ar[d]  \ar[r] & 0 \\ 
0 \ar[r] & W_n\sO_X \ar[r] \ar[d] & f_*W_n\sO_Y \ar[r] \ar[d] & \sQ_n \ar[r] \ar[d] & 0. \\
&  0  \ar[r] & \Ker(V^n : R^1 f_* W \sO_Y \to R^1 f_* W \sO_Y ) \ar[d] \ar[r] & \ar[d] \sE \ar[r] & 0 \\
&  & 0 & 0
}
\end{equation*}
As, we already know that $\sE$ and $\sQ_0$ are $p^N$ torsion, we obtain that $\sQ_n$ is indeed $p^{2N}$ torsion for every integer $n \geq 1$.
That is, we obtain the statement of \autoref{prop:wo-rationality_finite_level_it2}, for $M = 2N$.
\end{proof}

\noindent As a corollary of \autoref{prop:wo-rationality_finite_level} and \cite[Proposition 4.4.9]{Rulling_Chatzistamatiou_Hodge_Witt_and_Witt_rational}, we obtain:

\begin{proposition}\label{prop:properties_witt_rational}
Let $X$ be a variety over $k$.  Then the following assertions hold:
\begin{enumerate}
\item\label{it1:properties_witt_rational} the notion of $W\sO$-rationality is local in the \'etale topology,
\item\label{it2:properties_witt_rational} if $X' \to X$ is a universal homeomorphism then $X$ is W$\sO$-rational if and only if $X'$ is.
\end{enumerate}
\end{proposition}

\begin{proof}
The result \autoref{it2:properties_witt_rational} is proven \cite[Proposition 4.4.9]{Rulling_Chatzistamatiou_Hodge_Witt_and_Witt_rational}.  For \autoref{it1:properties_witt_rational}, we observe that \cite[I, Proposition 1.5.8]{Illusie_Complexe_de_de_Rham_Witt_et_cohomologie_cristalline} implies that for every surjective \'etale morphism $U \to X$ the associated map $W_n(U) \to W_n(X)$ is also surjective \'etale, and consequently for a quasi-resolution $f \colon Y \to X$ we obtain a cartesian diagram with surjective \'etale rows:
\[
\xymatrix{
    W_n Y \ar[d] & W_n(Y \times_X U) \ar[d]\ar[l] \\
    W_n X & W_n(U), \ar[l]^{W_n(\pi)}
}
\]
and $f_U \colon Y \times_X U \to U$ is a quasi-resolution according to \autoref{lem:etale_base_change_for_quasi_resolutions}.  Consequently, using flat base change (note that as above $W_n(\pi)$ is \'etale and hence flat), we see that
\[
R^if_{U,*} W_n\sO_{Y \times_X U} = W_n(\pi)^* R^if_* W_n\sO_Y.
\]
Since $W_n(\pi)$ is in fact faithfully flat, this implies that for every $M>0$ the integer $p^M$ annihilates $R^if_* W_n\sO_Y$ if and only if it annihilates $R^if_{U,*} W_n\sO_{Y \times_X U}$, which finishes the proof by \autoref{prop:wo-rationality_finite_level}.
\end{proof}

\begin{proposition}
\label{prop:separable_base_change_for_WO_rationality}
If $X$ is a W$\sO$-rational variety over $k$, and $K/k$ is a  separable (possibly infinite) algebraic field extenstion, then the variety $X_K$ is W$\sO$-rational.
\end{proposition}

\begin{proof}
Let $f \colon Y \to X$ be a W$\sO$-rational quasi-resolution, and 
let $\{K_\lambda\}_{\lambda \in I}$ be an ascending system of finite separable extension $k \subset K_\lambda \subset K$ such that $\bigcup_\lambda K_\lambda = K$. Set $\pi \colon X_K \to X$, $\psi \colon Y_K \to Y$, $\pi_{K_\lambda} \colon X_{K_\lambda} \to X$ and $\psi_{K_\lambda} \colon Y_{K_\lambda} \to Y$ to be the natural maps.   We shall prove that the base change $f_K \colon Y_K \to X_K$ is a W$\sO$-rational quasi-resolution.  It is a quasi-resolution by the arguments of  \autoref{lem:etale_base_change_for_quasi_resolutions}.  In order to see that it is also W$\sO$-rational we reason as follows.  Using \autoref{prop:wo-rationality_finite_level}, we choose an $M$ such that $p^M \cdot R^if_*W_n\sO_Y = 0$.  We claim it is also true for $f_K$ which concludes the proof using \autoref{prop:wo-rationality_finite_level} once again.  We consider the following diagrams, where by abuse of notation we identify every morphism $g$ with $W_n(g)$:
\[
\xymatrix{
    W_nY \ar[dd]_{f} & & W_n Y_K \ar[ll]_{\psi_K}\ar[dd]^{f_K} \ar[ld] \\
    & W_n Y_{K_\lambda} \ar[lu]^{\psi_{K_\lambda}} \ar[dd]|\hole \\
    W_nX & & W_n X_K \ar[ll]_(0.4){\pi_K} \ar[ld] \\
    &  W_n X_{K_\lambda} \ar[lu]^{\pi_{K_\lambda}} .\\
}
\]
Since the map $\pi_K$ is affine, we see that $p^M \cdot R^i f_{K,*} W_n\sO_{X_K} = 0$ if and only if $p^M \cdot \pi_{K,*} R^i f_{K,*} W_n\sO_{X_K} = 0$.  By the commutativity of the diagram above and the fact that $\psi_K$ is also affine, this is equivalent to the vanishing $p^M \cdot R^i f_*(\psi_{K,*} W_n \sO_{Y_K})$.  Since the forgetful functor from sheaves of rings to sheaves of sets preserves direct limits and $W_n$ is the $n$-fold self-product set theoretically, the pushforward $\psi_{K,*} W_n \sO_{Y_K}$ is in fact isomorphic to \[
\varinjlim_{\lambda} \psi_{K_\lambda,*}W_n \sO_{Y_{K_\lambda}}.
\] 
Now, since direct images commute with direct limits for quasi-compact quasi-separated morphisms (see \stacksproj{07U6}) we see that 
\begin{align*}
p^M \cdot R^i f_*(\psi_{K,*} W_n \sO_{Y_K}) & \isom p^M \cdot R^i f_*(\varinjlim_{\lambda} \psi_{K_\lambda,*} W_n \sO_{Y_{K_\lambda}})  \\
& \isom p^M \cdot \varinjlim_{\lambda} \pi_{K_\lambda,*} R^i f_{K_\lambda,*} W_n\sO_{Y_{K_\lambda}} = \varinjlim_{\lambda} p^M \cdot \pi_{K_\lambda,*} R^i f_{K_\lambda,*} W_n\sO_{Y_{K_\lambda}}.
\end{align*}
However, since the morphisms $X_{K_\lambda} \to X$ are in fact \'etale, by the argument in \autoref{prop:properties_witt_rational} all the terms in the last direct limit vanish which finishes the proof.
\end{proof}

\begin{remark}
We note that one would hope that a similar proof actually works for an arbitrary extension.  Unfortunately, a purely inseparable extension of height $k$ leads to the substitution of the exponent $M$ by $M+k$.  The limiting procedure cannot be therefore performed. 
\end{remark}


\subsection{Primer on $p$-adic cohomology}

For the sake of completeness, in this section we present the basic results concerning $p$-adic cohomology and direct images which are going to be necessary in our argument.  Most the results are contained in \cite[Exposé V ``Systemes projectifs $J$-adiques'']{SGA5} and \cite[Exposé VI ``Cohomologie $\ell$-adique'']{SGA5} by Jean-Pierre Jouanolou where the formalism of $\ell$-adic sheaves is developed.  We emphasize that the standard assumption $\ell \neq p$ showing up in $\ell$-adic cohomology considerations is not necessary for neither of the results in \emph{loc. cit}.

We begin with a short description of the necessary notion of Artin-Rees category $\AR(\ZZ_p)$ of $p$-adic sheaves.  We refer to \cite[Exposé V, \S 2.2]{SGA5} for the original account.  All the sheaves considered here are in  \'etale topology. The objects of the category are projective systems 
\[
\cF_\bullet = \left(\ \cF_n \ ,\ \phi^\cF_n \colon \cF_n \to \cF_{n-1}\ \right)_{n \in \ZZ}
\] of constructible (and hence torsion) sheaves of $\ZZ_p$-modules such that $\cF_n = 0$, for $n \ll 0$.  The morphisms are defined in a slightly subtle way, by the formula:
\[
\Hom_{\AR(\ZZ_p)}(\cF_\bullet,\cG_\bullet) = \varinjlim_{d \in \ZZ} \  \Hom_{\rm systems}(\cF_\bullet[d],\cG_\bullet)
\]
where $\cF[d]$ denotes the system defined by the association $\cF_\bullet[d]_n = \cF_{n+d}$, and the colimit is taken along the maps $\Hom_{\rm systems}(\cF_\bullet[d],\cG_\bullet) \to \Hom_{\rm systems}(\cF_\bullet[d+1],\cG_\bullet)$
induced by the diagram:
\[
\xymatrix{
    & & \cF_{n+d+1} \ar[r]^{\phi^{\cF}_{n+d+1}}\ar[dd]_{\phi^{\cF}_{n+d+1}} & \cF_{n+d} \ar[r]^{\psi_n}\ar[dd]_{\phi^{\cF}_{n+d}} & \cG_n \ar[dd]^{\phi^{\cG}_{n}} \\ 
    \psi = \Big(\psi_n \colon \cF_{n+d} \to \cG_n\Big)_{n \in \ZZ} \mapsto \Big(\psi_n \circ \phi^{\cF}_{n+d+1} \Big)_{n \in \ZZ} & & \\
    & & \cF_{n+d} \ar[r]_{\phi^{\cF}_{n+d}} & \cF_{n+d-1} \ar[r]_{\psi_{n-1}} & \cG_{n-1}.
}
\]
\noindent The right square in the diagram is in fact a morphism $\cF[d+1] \to \cF[d]$ of systems, and the above map is a composition therewith.  

We remark that the above definition is in fact equivalent to the localization of the category of systems of sheaves of $\ZZ_p$-modules satisfying the condition $\cF_n = 0$, for $n \ll 0$, by the subcategory of \emph{null systems}, that is, systems $\cN_{\bullet}$ for which the natural map $\cN[d] \to \cN$, induced by composition of $d$ transition maps, is zero for some $d \in \NN$ (see \cite[Exposé V, Proposition 2.4.3]{SGA5}).
\noindent A system $\cF_\bullet = (\cF_n)$ is \emph{strict} if $\cF_n = 0$, for $n<0$, and for $n \geq 0$ the sheaf $\cF_n$ is annihilated by $p^{n+1}$ and satisfies $\cF_n \otimes_{\ZZ/p^{n+1}} \ZZ/p^n \isom \cF_{n-1}$.  A strict system $\cF_\bullet$ is \emph{constructible} if all the sheaves $\cF_n$ are moreover constructible.

\begin{definition}[{$p$-adic sheaves}]
\label{defin:padic_sheaves}
We define the categories of $\ZZ_p$ and $\QQ_p$ sheaves as follows:
\begin{enumerate}
    \item the category of constructible $\ZZ_p$-sheaves is defined as the full subcategory of the category $\AR(\ZZ_p)$ of systems isomorphic to a strict constructible system.
    \item the category of constructible \emph{$\QQ_p$-sheaves} (see \cite[Exposé VI, Définition 1.4.3]{SGA5}) is defined as the localization of the category of $\ZZ_p$-sheaves by the Serre subcategory of sheaves annihilated by a fixed power of $p$ (see \autoref{remark:localization} for a brief description of the localization procedure). 
\end{enumerate}
Let $f \colon X \to Y$ be a proper finite type morphism of locally noetherian schemes.  The most important reason for the introduction of the above admittedly complex definition is that it allows for the direct images of a $\ZZ_p$-sheaf $\cF = (\cF_n)_{n \in \NN}$ to be well-defined using the formula:
\[
R^if_*\cF = (R^if_*\cF_n)_{n \in \NN}.
\]
We remark that the above higher direct image systems are often not strict constructible themselves.  However, by the proper base change for $\ZZ/p^n$-constructible sheaves and the main result of \cite[Exposé VI, \S 2.2]{SGA5}, they are Artin--Rees isomorphic to a strict constructible system.  In fact, the only condition which is likely not satisfied is strictness, which however holds after an Artin--Rees isomorphism.  Moreover, if $f \colon \overline{X} \to \overline{k}$ is a structure morphism of a base change over the algebraic closure of a variety $X/k$, we recover the usual definition of $p$-adic cohomology using \cite[Exposé VI, 1.2 ``$\ZZ_\ell$ faisceaux constants tordus constructible'']{SGA5}.  We note that in this case all sheaves $R^if_*\ZZ/p^n$ are in fact torsion constant and hence clearly satisfy requirements of \emph{loc. cit}.

\end{definition}

We shall need the following standard results concerning behaviour of constructible $\ZZ_p$-sheaves.  We remark that although the following results are stated in \cite[Exposé VI]{SGA5} only for $\ZZ_p$-sheaves, the author explains in \cite[Exposé VI, 2.2.5]{SGA5} that all the results can be carried through to the $\QQ_p$ setting using straightforward localization.

\begin{proposition}[{\cite[Exposé VI, Proposition 2.2.4]{SGA5}}]
\label{prop:leray_spectral_seq_QQ_p}
Let $f \colon X \to Y$ and $g \colon Y \to Z$ be proper finite type morphisms of locally noetherian schemes.  Then for every constructible $p$-adic sheaf $\cF$ there exists a convergent Leray spectral sequence: 
\[
E^{pq}_2 = R^pg_*(R^qf_*\cF)\Rightarrow R^{p+q}(g \circ f)_* \cF.
\]
In particular, if $Y$ is a proper variety defined over an algebraically closed field and $f \colon X \to Y$ is proper finite type morphism then there exists a convergent Leray spectral sequence:
\[
E^{pq}_2 = H^p_{\et}(Y,R^qf_*\QQ_p) \Rightarrow H^{p+q}_{\et}(Y,\QQ_p).
\]
\end{proposition}

\noindent In fact, the result in \emph{loc. cit} works for higher direct images with compact support but we will only need it in the stated generality.

\begin{proposition}[{\cite[Exposé VI, 2.2.3 B)]{SGA5}}]
\label{prop:QQ_p_proper_base_change}
Suppose that the diagram
\[
\xymatrix{
    X \ar[d]_f & X' \ar[l]_v \ar[d]^g \\
    Y & \ar[l] Y' \ar[l]^u
}
\]
is a cartesian square of locally noetherian schemes.  Assume that the morphism $f$ is proper and of finite type.  Then for every constructible $\ZZ_p$-sheaf $\cF$ (resp. $\QQ_p$-sheaf) the natural base change map:
\[
u^*R^if_*\cF \to R^ig_*(v^*\cF)
\]
is an isomorphism.
\end{proposition}

\begin{theorem}[{\cite[Théorème 5.1]{Chambert_Loir_Points_Rationnels_et_Groupes_Fondamentaux}}]
\label{thm:clr_rigid_zero_padic}
For every proper finite type scheme $X$ defined over an algebraically closed field the following identification holds true:
\[
H^\bullet_{\et}(X,\QQ_p) \otimes K \isom H^\bullet_{\rm rig}(X)^{=0},
\]
where $K:=\Frac(W(k))$.
\end{theorem}

\noindent As a direct combination of \autoref{thm:bbe_rigid_witt_vector} and \autoref{thm:clr_rigid_zero_padic} we obtain the following:

\begin{corollary}
\label{cor:p_adic_slope_zero_of_witt_vector}
For every finite type scheme $X$ defined over an algebraically closed field the following identification holds true:
\[
H^\bullet_{\et}(X,\QQ_p) \otimes K \isom H^\bullet(X,W\sO_{X,\QQ})^{=0},
\]
where $K:=\Frac(W(k))$.
\end{corollary}

\subsection{$p$-adic cohomology and W$\sO$-rational morphisms}

In \autoref{cor:p_adic_slope_zero_of_witt_vector} we saw that the absolute $p$-adic cohomology is closely related to the Witt vector cohomology.  The subsequent lemma provides a similar relation in the relative setting.  We only present this simplistic version of the theorem, because a fully-fledged comparison will most likely require application of rigid cohomology and $F$-isocrystals (as do the full versions of \autoref{thm:bbe_rigid_witt_vector} and \autoref{thm:clr_rigid_zero_padic}).  

The main tool that we utilize is the system of Witt Artin--Schreier sequences exact in \'etale topology:
\begin{align}
\label{eq:witt_artin_schreier}
0 \ra \ZZ/p^n \ra W_n\sO_X \xrightarrow{1 - F} W_n\sO_X \ra 0.
\end{align}

\begin{lemma}
\label{lemma:wo_rational_implies_QQp_rational}
Let $f \colon X \to Y$ be a proper W$\sO$-rational morphism (see \autoref{remark:some_remarks_on_WO_rationality}).  Then $f_*\QQ_p = \QQ_p$ and $R^if_*\QQ_p = 0$ for every $i>0$.
\end{lemma}

\begin{proof}
Looking at degree $i>1$ terms of the long exact sequences of cohomology associated to \eqref{eq:witt_artin_schreier} we obtain exact sequences:
\begin{equation}
\label{eq:wo_rational_implies_QQp_rational:LES}
\cdots \ra R^{i-1}f_*W_n\sO_X \ra R^if_*\ZZ/p^n \ra R^if_*W_n\sO_X \ra \cdots.
\end{equation}
Since $f$ is a W$\sO$-rational, according to \autoref{prop:wo-rationality_finite_level}, the left and right terms of \autoref{eq:wo_rational_implies_QQp_rational:LES} are annihilated by a fixed number $p^M$ and hence for every $n>0$ the middle term is killed by $p^{2M}$.  Using \autoref{defin:padic_sheaves}, this directly implies that $R^if_*\QQ_p = 0$, for $i>1$. For low degrees, for every $n>0$ we consider the diagram of sheaves on the \'etale site of $Y$:
\begin{equation}
\label{eq:wo_rational_implies_QQp_rational:big}
\xymatrix{
    & 0 \ar[d] & 0 \ar[d] & & 0 \ar[d] \\
    0 \ar[r] & (\ZZ/p^n)_Y \ar[r]\ar[d] & W_n\sO_Y \ar[rr]^{1 - F_Y}\ar[d] & & W_n\sO_Y \ar[r]\ar[d] & 0 \ar[d] \\
    0 \ar[r] & f_*(\ZZ/p^n)_X \ar[r]\ar[d] & f_*W_n\sO_X \ar[rr]^{f_*(1 - F_X)}\ar[d] & & f_*W_n\sO_X \ar[r]\ar[d] & \coker\left(f_*(1 - F_X)\right) \ar[d] \ar[r] & 0  \\
    0 \ar[r] & K_n \ar[d] \ar[r] & C_n \ar[rr]\ar[d] & & C_n \ar[d]\ar[r] & C_n/f_*(1-F_X) C_n \ar[r] & 0 \\
    &  0 & 0 & & 0
} 
\end{equation}
A priori, we only know that the leftmost three non-zero columns, as well as the middle row and the row above that are exact in \autoref{eq:wo_rational_implies_QQp_rational:big}. Then, Snake lemma applied to the first three non-zero entries of these two rows yields, that the row below the middle row is also exact. Finally, Snake lemma applied to the second and the third non-zero column yields that also the rightmost column is exact. So, eventually we obtain that all columns and rows of \autoref{eq:wo_rational_implies_QQp_rational:big} are exact.

By the middle columns of \autoref{eq:wo_rational_implies_QQp_rational:big} we obtain that $C_n$ are annihilated by $p^{M}$ for every $n>0$. Then, the lowest non-zero row of \autoref{eq:wo_rational_implies_QQp_rational:big} implies that $K_n$ are annihilated by $p^{2M}$ for every $n>0$. Again, using \autoref{defin:padic_sheaves} and the fact that localization by Serre subcategories is an exact functor, this implies that $f_* \QQ_p = \QQ_p$.  Moreover, for $i=1$, we observe that \autoref{eq:wo_rational_implies_QQp_rational:LES} and  the middle row of \autoref{eq:wo_rational_implies_QQp_rational:big}  lead to the exact sequence
\[
0 \ra \coker\left(f_*(1 - F_X)\right) \ra R^1f_*(\ZZ/p^n)_X \ra R^1f_*W_n\sO_X. 
\]
Since both $\coker\left(f_*(1 - F_X)\right)$, which injects into $C_n/f_*(1-F_X) C_n$, and $R^1f_*W_n\sO_X$ are $p^{M}$-torsion, the sheaf $R^1f_*(\ZZ/p^n)_X$ is $p^{2M}$-torsion for every $n>0$ and hence $R^1f_*\QQ_p = 0$.  This finishes our proof.
\end{proof}

\noindent In the spirit of \autoref{remark:some_remarks_on_WO_rationality}, the morphism satisfying the properties in stated in \autoref{lemma:wo_rational_implies_QQp_rational} is called \emph{$\QQ_p$-rational}. That is a proper $f : X \to Y$ is $\QQ_p$-rational if $f_*\QQ_p = \QQ_p$ and $R^if_*\QQ_p = 0$ for every $i>0$.

\begin{proposition}
\label{prop:geom_generic_fibres_QQp_rational}
Let $f \colon X \to T$ be a flat morphism of varieties defined over $k$.  Assume that $X$ is W$\sO$-rational.  Then the geometric generic fibre $X_{\overline{\eta}}$ of $f$ admits a $\QQ_p$-rational quasi-resolution.
\end{proposition}

\begin{proof}
Let $\overline{\eta}$ be a geometric generic point of $T$.  Let $\xi \colon Y \to X_{\overline{\eta}}$ be a quasi-resolution.  As $Y$ and $\xi$ are of finite type over $k\left(\overline{\eta} \right)$, both are defined over a finite extension $L$ of $k(\eta)$. Taking the normalization $U_{\pre}$ of $T$ in $L$, and restricting to the open set where the separable part of $U_{\pre} \to T$ is \'etale we obtain:
\begin{itemize}
\item a morphism $U \to T$, which is a composition of an \'etale map and a universal homeomorphism, and
\item  a quasi-resolution $\pi \colon Z \to X_U$ such that $\xi = \pi_{\overline{\eta}}$.
\end{itemize}  
Using \autoref{prop:properties_witt_rational}, the scheme $X_U$ is W$\sO$-rational and hence $\pi$ is a W$\sO$-rational morphism.  By \autoref{lemma:wo_rational_implies_QQp_rational}, the morphism $\pi$ is also $\QQ_p$-rational.  Applying the proper base change (see \autoref{prop:QQ_p_proper_base_change}) this property is inherited by $\xi = \pi_{\overline{\eta}}$ which finishes the proof.
\end{proof}

\begin{remark}
Although \autoref{prop:geom_generic_fibres_QQp_rational} is  sufficient for our purpose,  it is  admittedly a workaround. A clearer claim would be that the geometric generic fibre is in fact W$\sO$-rational itself.  We were unfortunately unable to prove it.  The main obstacle is the lack of a proper base change type theorem, at the geometric generic point, for Witt vector higher direct images.
\end{remark}

%% file: top_Witt.tex
\section{Top Witt vector and $p$-adic cohomology}

In this part of the paper we present a few results concerning top Witt vector cohomology of varieties.  
We begin with the following simple technical result which allows us to describe the top Witt vector cohomology in terms of coherent cohomology under arithmetic assumption of ordinarity.  This result is most likely standard but we could not find a precise reference.

\begin{proposition}\label{lem:wo_o_cohomology_ordinary}
Let $X$ be a proper variety of dimension $d = \dim X$ defined over a perfect field $k$ of characteristic $p>0$.  Assume that  $X$ is weakly ordinary, that is, the natural Frobenius action on $H^d(X,\sO_X)$ is bijective, and that $H^{d-1}(X,\sO_X) = 0$.  Then $H^d(X,W\sO_{X,\QQ})$ is a $K$-vector space of dimension $g = \dim_k H^d(X,\sO_X)$, and 
\[
H^d_{\et}(X,\QQ_p) \otimes_{\QQ_p} K \isomto H^d(X,W\sO_{X,\QQ}),
\]
where $K:= \Frac W(k)$.
\end{proposition}

\begin{proof}
First, we observe that by a straightforward inductive argument $H^{d-1}(X,W_n\sO_X) = 0$, for every $n \geq 1$.  Then, we claim that the map $F \colon H^d(X,W_n\sO_X) \to H^d(X,W_n\sO_X)$ induced by the Witt vector Frobenius is an isomorphism.  This follows directly by induction on $n$ from the five lemma applied to following diagram, the rows of which are exact by the vanishing of finite length Witt vector cohomology in degree $d-1$. 
\[
\xymatrix{
0 \ar[r] & H^d(X,\sO_X) \ar[r]^(0.45){V} \ar[d]^F & H^d(X,W_{n+1}\sO_X) \ar[r]\ar[d]^F & H^d(X,W_n\sO_X) \ar[r]\ar[d]^F & 0 \\
0 \ar[r] & H^d(X,\sO_X) \ar[r]^(0.45){V} & H^d(X,W_{n+1}\sO_X) \ar[r] & H^d(X,W_n\sO_X) \ar[r] &  0,
}
\]
 Taking limit with respect to $n$ this implies that the Frobenius endomorphism of $H^d(X,W\sO_X)$ is also an isomorphism.  Using the vanishing of $H^{d-1}(X,\sO_X)$ again, we see that connecting homomorphism $H^{d-1}(X,\sO_X) \to H^d(X,W\sO_X)$ in the long exact sequence associated with
\[
0 \ra W\sO_X \xrightarrow{V} W\sO_X \ra \sO_X \ra 0.
\]
is zero.  Therefore the map $V \colon H^d(X,W\sO_X) \to H^d(X,W\sO_X)$ is injective and hence $p = FV$ is injective as well.  By \autoref{prop:witt_vectors_and_lim}, the natural map $H^d(X,W\sO_X) \to \varprojlim_{n} H^d(X,W_n\sO_X)$ is an isomorphism, which implies that $\bigcap_n V^n H^d(X,W\sO_X) = \{0\}$ and therefore $H^d(X,W\sO_X)$ is $V$-adically and hence $p$-adically separated.  Consequently, using Nakayama's lemma for the complete ring $W(k)$ (see \cite[Theorem 8.4]{Matsumura_Commuatative_ring_theory}), we see that the torsion-free and $p$-adically separated module $H^d(X,W\sO_X)$ is in fact free of rank:
\[
\dim_k H^d(X,W\sO_X)/pH^d(X,W\sO_X) = \dim_k H^d(X,W\sO_X)/VH^d(X,W\sO_X) = \dim_k H^d(X,\sO_X),
\]
as desired.  After localization, this yields the result concerning Witt vector cohomology.  

For the claim related to $p$-adic cohomology, we take a limit of the cohomology exact sequences associated to \eqref{eq:witt_artin_schreier} in order to obtain a short exact sequence (note that we use the vanishing in degree $d-1$):
\[
0 \ra H^d_{\et}(X,\ZZ_p) \ra H^d(X,W\sO_X) \xrightarrow{1 - F} H^d(X,W\sO_X) \ra 0,
\]
which means that $H^d_{\et}(X,\ZZ_p)$ are in fact Frobenius stable elements in the $F$-crystal $H^d(X,W\sO_X)$.
Since, as proven above, the morphism $F$ acting on $H^d(X,W\sO_X)$ is a $p$-linear automorphism, we may now use the standard result on $F$-crystals (see \cite[2.1.2]{Katz_Travaux_de_Dwork}) which implies that $H^d(X,W\sO_X)$ is generated by Frobenius stable elements.  This yields the claim after localization.
\end{proof}

We now proceed to the proof that both Witt vector and $p$-adic cohomology vanishes for uniruled W$\sO$-rational varieties.  We precede the actual statement with a few lemmata.  Before stating them we recall that: if $X$ and $Y$ are integral, and $Y$ is normal, then $u \colon X \to Y$ is a universal homeomorphism if and only if it is finite, surjective and purely inseparable (see \cite[Prop. (3.5.8)]{EGAI}).

\begin{lemma}[{\cite[Lemma 4.1.6 \& Lemma 4.2.4]{Rulling_Chatzistamatiou_Hodge_Witt_and_Witt_rational}}]
\label{lemma:witt_cohomology_universal_homeo}
Let $f \colon Y \to X$ be a morphism of normal varieties over $k$.  Then the following statements hold true for every $i>0$:
\begin{itemize}
    \item if $f$ is finite then $f^* \colon H^i(Y,W\sO_{Y,\QQ}) \to H^i(X,W\sO_{X,\QQ})$ is injective,
    \item if $f$ is a universal homeomorphism  then $f^* \colon H^i(Y,W\sO_{Y,\QQ}) \to H^i(X,W\sO_{X,\QQ})$ is an isomorphism.
\end{itemize}
\end{lemma}

\begin{lemma}\label{lemma:support_witt_cohomology}
Let $f \colon X \to Y$ be a birational morphism of varieties of dimension $d$.  Then for every $i>0$ we have 
\[
\dim {\rm Supp}(R^if_*W\sO_{X,\QQ}) < d-i.
\]
\end{lemma}

\begin{proof}
We first recall the standard proof that $\dim {\rm Supp}(R^if_*\sO_X) < d-i$.  Let $\eta \in Y$ be the generic point of ${\rm Supp}(R^if_*\sO_X)$.  Localizing and using the formal function theorem, we see that 
\[
(R^if_*\sO_X)^{\wedge}_{\eta} = \lim_{n \to 
 \infty} H^i(X_n,\sO_{X_n}),
\]
for $X_n = Y \times_X \Spec(\sO_{X,\eta}/\mathfrak{m}^{n+1}_{\eta})$ where $\mathfrak{m}_\eta \subset \sO_{X,\eta}$ is the maximal.  This implies that $i \leq \dim X_n = \dim X_0$.  Moreover, since $\eta$ is contained in the image of the exceptional set of $f$ we see that $\dim {\rm Supp}(R^if_*\sO_X) + \dim X_0 < d$.  Combining those inequalities we obtain the claim.

We now proceed to the proof concerning $(R^if_*W\sO_{X,\QQ})$.  Set $S_i = {\rm Supp} (R^if_*\sO_X)$  We begin by observing that a simple inductive argument with respect to the parameter $n$ coming from the exact sequence:
\[
\cdots \ra R^if_*\sO_X \xrightarrow{V^n} R^if_*W_{n+1}\sO_X \ra R^if_*W_n\sO_X \ra \cdots,
\]
implies that ${\rm Supp} (R^if_*W_n\sO_X) \subset S_i$.  Moreover, from the exact sequence
\[
\cdots \ra R^{i-1}f_*W_{n+1}\sO_X \xrightarrow{V} R^{i-1}f_*W_n\sO_X \ra R^if_*\sO_X \ra \cdots,
\]
we see that the maps $R^{i-1}f_*W_{n+1}\sO_X \ra R^{i-1}f_*W_n\sO_X$ are surjective on $X \setminus S_i$.  Consequently, the associated projective system satisfies Mittag-Leffler condition along $X \setminus S_i$ and hence 
\[
{\rm Supp}(R^1\lim R^{i-1}f_*W_n\sO_X) \subset S_i.
\]  
Combination of the above two observations with \autoref{prop:witt_vectors_and_lim} finishes the proof.
\end{proof}

\begin{definition}
We say that a morphism of normal varieties $X \to Y$ is \emph{quasi-birational} if it is surjective purely inseparable and generically finite.  Alternatively, taking Stein factorization, it is a composition of a universal homeomorphism and a birational morphism.
\end{definition}

\begin{lemma}
\label{proposition:quasi_birational}
Let $f \colon X' \to X$ be a quasi-birational morphism of normal varieties of dimension $d$.  Then the map $f^* \colon H^d(X,W\sO_{X,\QQ}) \to H^d(X',W\sO_{X',\QQ})$ is surjective.
\end{lemma}

\begin{proof}
First, taking Stein factorization and using \autoref{lemma:witt_cohomology_universal_homeo}, we see that we may assume that $f$ is birational.  In this situation, we conclude by Leray spectral sequence and \autoref{lemma:support_witt_cohomology}.  More precisely, we observe that the top diagonal of the $E_2$ page of the spectral sequence consists of a single entry $H^d(X,W\sO_{X,\QQ})$, which clearly yields the result.
\end{proof}

\begin{proposition}
\label{lem:uniruled_witt_vector}
Let $X$ be a W$\sO$-rational normal uniruled variety of dimension $d = \dim X$ defined over an algebraically closed field of characteristic $p>0$.  Then $H^d(X,W\sO_{X,\QQ}) = 0$ and $H^d_{\et}(X,\QQ_p) = 0$
\end{proposition}

\begin{proof}
Since $X$ is uniruled there exists a generically finite rational map $g \colon Z \times \PP^1 \dashrightarrow X$.  Taking the quasi-resolution of the closure of the graph $\Gamma_g \subset Z \times \PP^1 \times X$ we obtain a proper generically finite map $Y \to X$ with $H^d(Y,W\sO_{Y,\QQ}) = 0$ by \autoref{proposition:quasi_birational}. 

\[
\xymatrix{
  Z \times \PP^1 \ar@{-->}[d] & Y \ar[l]\ar[d] & Y' \ar[l]\ar[d] \\
  X \ar@{=}[r] & X & \ar[l] X'
}
\]

\noindent Taking the flattening of the morphism $Y \to X$, performing base change along a quasi-resolution of the base and potentially taking the normalization of the dominant component, we obtain a proper finite morphism $Y' \to X'$ with $Y'$ normal and quasi-birational to $Y$ and $X'$ a quasi-resolution of $X$.  By applying \autoref{proposition:quasi_birational} again, we see that $H^d(Y',W\sO_{Y',\QQ}) = 0$, which implies that $H^d(X',W\sO_{X',\QQ}) = 0$ using \autoref{lemma:witt_cohomology_universal_homeo}.  Since $X$ is W$\sO$-rational, this means that $H^d(X,W\sO_{X,\QQ}) = H^d(X',W\sO_{X',\QQ}) = 0$, and hence the proof concerning W$\sO$-cohomology is finished.  The claim about $p$-adic cohomology is now a direct consequence of \autoref{cor:p_adic_slope_zero_of_witt_vector}.
\end{proof}

%% file: CY.tex

\section{Weakly ordinary Calabi--Yau varieties}

In the following section, we prove our main result stating that weakly ordinary varieties with trivial canonical class are not uniruled.  We begin with a few preliminary results.


\subsection{Maximal rationally chain connected fibrations}\label{ss:mrcc_fibrations}

We first recall some basic results concerning maximal rationally chain connected (MRCC for short) fibrations.  All the results are well-explained in \cite[IV.5]{Kollar_Rational_curves_on_algebraic_varieties}.   

\begin{definition}[Rationally chain connected fibration]
Suppose $X$ is a normal scheme over a field $k$.  Let $X^{\circ} \subseteq X$ be an open subset, and let $\pi^{\circ} \colon X^{\circ} \to S^{\circ}$ be a morphism.
\begin{enumerate}
	\item We say that $f^{\circ}$ is a \emph{rationally chain connected fibration} if it is a proper morphism satisfying $f^{\circ}_*\cO_{X^{\circ}} = \cO_{S^{\circ}}$ such that all the geometric fibres are rationally chain connected.
	\item We say that $f^{\circ}$ is a \emph{maximal rationally chain connected} if for every other open $X' \subseteq X$ and a rationally chain connected fibration $f' \colon X' \to S'$ there exists a rational map $\pi \colon S' \dasharrow S^{\circ}$ such that $\pi \circ f' = f^{\circ}$. 
\end{enumerate} 
\end{definition}

\begin{theorem}[{\cite[Theorem IV.5.2 \& Complement IV.5.2.1]{Kollar_Rational_curves_on_algebraic_varieties}}]\label{thm:existence_mrc_fibrations}
Let $X$ be a normal proper variety.  Then a maximal rationally chain connected fibration $f^{\circ} \colon X^{\circ} \to S^{\circ}$ exists.  Moreover, the fibration of a uniruled variety $X$ is non-trivial.
\end{theorem}

\begin{proof}
The first claim is proven in \emph{loc. cit}.  In \cite[Proposition A.4]{Tanaka_Behavior_of_canonical_divisors_under_purely_inseparable_base_changes} Tanaka reproves the second claim also stating that it is in fact implicit in \cite[Theorem IV.5.2]{Kollar_Rational_curves_on_algebraic_varieties}.
\end{proof}

\subsection{\'Etale fundamental groups}

%
%
%

Here we present a few results concerning fundamental groups.  We begin with the recollection of the standard corollary of the Artin--Schreier sequence.

\begin{lemma}\label{lem:fundamental_group_artin_schreier}
Let $X$ be a proper globally $F$-split variety defined over an algebraically closed field of characteristic $p>0$ such that $H^1(X,\cO_X) \neq 0$.  Then $H^1(X,\FF_p) \neq 0$, and consequently $X$ admits a non-trivial $\FF_p$-covering.
\end{lemma}
\begin{proof}
Applying Artin-Schreier sequence of \'etale sheaves:
\[
0 \ra \FF_p \ra \cO_X \xrightarrow{F - \id} \cO_X \ra 0,
\]  
along with the \'etale descent, we obtain the long exact sequence of cohomology
\[
\xymatrix{
  \cdots \ar[r] & H^1_{\et}(X,\FF_p) \ar[r] & H^1(X,\cO_X) \ar[r]^{F - \id} & H^1(X,\cO_X) \ar[r] & \cdots,
}
\] 
where $F \colon H^1(X,\cO_X) \to H^1(X,\cO_X)$ is the natural $p$-linear map induced by the Frobenius on $X$.  Using standard results in semi-linear algebra, we take a decomposition:
\[
H^1(X,\cO_X) = H^1(X,\cO_X)^{\rm ss} \oplus H^1(X,\cO_X)^{\rm nil},
\]
into a semistable part $H^1(X,\cO_X)^{\rm ss}$ generated by $F$-stable elements, and nilpotent part \linebreak $H^1(X,\cO_X)^{\rm nil}$ where $F$ is nilpotent.  Since $X$ is globally $F$-split, the map $F$ is is bijective, and hence $H^1(X,\cO_X)^{\rm ss} = H^1(X,\cO_X) \neq 0$.  Consequently, we see that $H^1_{\et}(X,\FF_p) \neq 0$, which implies that $X$ admits a non-trivial $\FF_p$-covering.
\end{proof}



\begin{theorem}[{\cite[Théorème]{Chambert_Loir_groupe_fondamental_rationnellement_connexes}}]
\label{thm:fundamental_group_rcc}
Let $X$ be a proper normal rationally chain connected variety defined over an algebraically closed field of characteristic $p>0$.  Then $\piet(X)$ is a finite group.
\end{theorem}

\noindent Even though the statement in \emph{loc. cit.} actually concerns varieties over algebraically closed fields, taking a perfect closure does not change the fundamental group.


\begin{remark}
In the reference applied in \autoref{thm:fundamental_group_rcc}, the author also gives a short indication that some of the techniques included in his paper \cite{Chambert_Loir_Points_Rationnels_et_Groupes_Fondamentaux} could lead to a slightly more detailed result, potentially useful for our purposes, that the order of the fundamental group is in fact prime to $p$.  We were not able to fully verify the claim so we only apply the above result, which is in the end sufficient for our main argument. 
\end{remark}

\subsection{Main result}

We now approach the main theorem of the present paper. 

\begin{theorem}\label{thm:weakly_ordinary_cy_not_uniruled}
Let $X$ be a normal globally $F$-split and W$\cO$-rational proper variety defined over $k$.  
Suppose that either:
\begin{enumerate}
     \item\label{thm:weakly_ordinary_cy_not_uniruled:it1} $X$ satisfies Serre's condition $S_3$ and $K_X \sim 0$, or
     \item\label{thm:weakly_ordinary_cy_not_uniruled:it2} $X$ is smooth and $K_X$ is numerically trivial.
\end{enumerate} 
Then $X$ is not geometrically uniruled.
\end{theorem}

\begin{proof}
We first make a straightforward reduction of \autoref{thm:weakly_ordinary_cy_not_uniruled:it2} to \autoref{thm:weakly_ordinary_cy_not_uniruled:it1}.  Indeed, by the standard bijection between $F$-splittings and appropriate sections of $(1-p)K_X$ (see, e.g., \cite[4.2]{Schwede_Smith_Globally_F_regular_an_log_Fano_varieties}), we have a linear equivalence $(1-p)K_X \sim 0$ since $K_X$ is numerically trivial.  Consequently, there exists an \'etale $\ZZ/(p-1)$-cyclic covering $\phi: X' \to X$ trivializing $K_X$ which satisfies all the necessary properties: one only needs to show that $X'$ is globally $F$-split, for which the splitting of $\sO_X' \to F_* \sO_{X'}$ is given by the pullback via $\phi$ of the splitting of $\sO_X \to F_* \sO_X$.  

We now proceed to the proof of \autoref{thm:weakly_ordinary_cy_not_uniruled:it1}.  We suppose that $X$ is a geometrically uniruled variety of dimension $d$, and we seek a contradiction.  Taking a separable extension of the base field we may assume that $k$ is algebraically closed, $X$ is uniruled and all the assumed properties are preserved.  Indeed, we use the following references: \stacksproj{0C3M} for normality, \autoref{prop:separable_base_change_for_WO_rationality} for W$\cO$-rationality, \cite[Corollary 2.5]{Gongyo_Li_Patakfalvi_Schwede_Tanaka_Zong_On_rational_connectedness_of_globally_F_regular_threefolds} for $F$-splitting and \cite[Corollaire 6.4.2]{Grothendieck_Elements_de_geometrie_algebrique_IV_II}  for $S_3$ property.  By \autoref{thm:existence_mrc_fibrations} there exists an open subset $X^{\circ} \subseteq X$ and a non-trivial MRCC fibration $f^{\circ} \colon X^{\circ} \to S^{\circ}$. 

We claim that \emph{all the stated properties of $X$ are inherited to the geometric generic fibre $X_{\overline{\eta}}$ of $f^{\circ}$, except possibly $W \sO$-rationality}. Indeed, by \autoref{cor:f_split_normal_fibres}, $X_{\overline{\eta}}$ is globally $F$-split and normal.  Then \autoref{prop:property_S_k_inherited_by_fibres} directly implies that the geometric generic fibre is $S_3$. Other properties follow from adjunction, which concludes the proof of our claim.

By the definition of the MRCC fibration, we know that $X_{\overline{\eta}}$ is rationally chain connected.
Furthermore, by \autoref{prop:geom_generic_fibres_QQp_rational}, $X_{\overline{\eta}}$ admits a $\QQ_p$-rational quasi-resolution.  Consequently, by replacing $X$ be $X_{\overline{\eta}}$, we may assume that:
\begin{enumerate}
    \item $X$ is normal and $S_3$,
    \item $X$ is globally $F$-split,
    \item $K_X$ is trivial,
    \item $X$ is rationally chain connected, and
    \item $X$ admits a $\QQ_p$-rational quasi-resolution.
\end{enumerate}

We now claim that we can also assume that $H^1(X,\cO_{X}) = 0$.  First, we observe that by \autoref{thm:fundamental_group_rcc} the \'etale fundamental group of $X$ is finite.  Using \autoref{lem:fundamental_group_artin_schreier}, if the necessary vanishing $H^1(X,\cO_{X}) = 0$ is not satisfied, there exists an \'etale $\ZZ/p$-covering $X'$ which satisfies all the assumptions concerning $X$ including rational chain connectedness and existence of a $\QQ_p$-rational quasi-resolution.  We note that existence of an appropriate quasi-resolution follow from  \autoref{lem:etale_base_change_for_quasi_resolutions} and \autoref{prop:QQ_p_proper_base_change}.  We now repeatedly substitute $X$ with $X'$, and observe that this process needs to terminate because the order of the \'etale fundamental group decreases along the way.

Equipped with the above, we now use \autoref{prop:duality_for_S_k} ($X$ is $S_3$) in order to see that $H^{d-1}(X,\cO_X)$ vanishes too.  Consequently, \autoref{lem:wo_o_cohomology_ordinary} gives non-vanishings $H^d_{\et}(X,\QQ_p) \neq 0$.
We now consider the $\QQ_p$-rational quasi-resolution $Y \to X$ whose existence is postulated above.  By $\QQ_p$-rationality and the Leray spectral sequence for $p$-adic cohomology (see \autoref{prop:leray_spectral_seq_QQ_p}) we know that
\[
H^d_{\et}(Y,\QQ_p) \isomto H^d_{\et}(X,\QQ_p) \neq 0.
\]
Since $Y \to X$ is a quasi-birational morphism, one easily sees that $Y$ is uniruled and consequently the non-vanishing above violates \autoref{lem:uniruled_witt_vector} yielding the desired contradiction.  The proof is thus finished.
\end{proof}

\noindent Using the results of Langer \cite[Corollary 3.3]{Langer_Generic_positivity_and_foliations_in_positive_characteristic}, we obtain the following corollary.

\begin{corollary}
Let $X$ be a weakly ordinary projective variety with trivial canonical class.  Then the tangent bundle $\cT_X$ is strongly semistable with respect to some (equiv. every) polarization.
\end{corollary}


\noindent We now provide three classes of examples justifying certain steps in the above reasoning.  

\begin{remark}
We attempted to prove the above theorem in the more general setting where $X$ is strongly $F$-regular and $K_X$ is only numerically trivial.  Under this assumptions the $\ZZ/(p-1)$-cyclic covering trivializing $K_X$ appearing in the proof is only quasi-\'etale.  Unfortunately, neither W$\cO$-rationality nor rational chain connectedness is preserved under quasi-\'etale maps without any additional assumptions.  For the statement concerning rational chain connectedness, one considers the following example.  Let $N \in \NN$ be coprime to $p$ and $G = \factor{\ZZ}{N \ZZ}$.  Let $C$ be a smooth projective curve of positive genus defined over an algebraically closed field and equipped with a $G$-action such that $C/G \isom \PP^1$.  Moreover, assume that $G$ acts linearly on $\PP^2$ via diagonal matrices ${\rm diag}(\xi^a,\xi^b,\xi^c)$, where $\xi$ is the $N$-th primitive root of unity and the numbers $a$, $b$, $c$ and $N$ are pairwise coprime.  The variety $\PP^2 \times C$ is clearly not rationally chain connected.  We claim that the quotient $(\PP^2 \times C)/ G$ is rationally chain connected and the map $\PP^2 \times C \to (\PP^2 \times C) / G$ is quasi-\'etale.  The second claim is clear, because the fixed points of the action on the product occur in codimension two.  For the first, we observe that the natural map 
\[
(\PP^2 \times C)/G \to C/G = \PP^1
\]
is generically a $\PP^2$-bundle over $\PP^1$ admitting a section by Tsen's theorem stating that the Brauer group of a function field of a curve over an algebraically closed field is trivial and hence all $\PP^n$-bundles admit a section (see \stacksproj{03RC}).  This example also shows that fundamental groups behave wildly under quasi-\'etale morphisms.  For example, we see that $\pi_1\left((\PP^2 \times C)/G\right)$ is finite, and $\pi_1(\PP^2 \times C) \isom \pi_1(C)$ is infinite highly non-commutative.

\end{remark}

\begin{example}\label{example:singular_enriques_surfaces}
One could hope that the non-vanishing of top rational Witt vector cohomology holds for every weakly ordinary Calabi--Yau variety.  Unfortunately, this is not the case.  We claim that the counterexample is provided by singular Enriques surfaces in characteristic two.  Indeed, from the Bombieri--Mumford classification (see \cite[Section 1]{Liedtke_Arithmetic_moduli_and_lifting_of_Enriques_surfaces} for the details) we know that for such a surface $S$ the canonical divisor is trivial and moreover the vector spaces $H^1(S,\cO_S)$ and $H^2(S,\cO_S)$ are one-dimensional with bijective Frobenius action.  This implies that such $S$ is weakly ordinary and Calabi--Yau.  However, in \cite[II, Proposition 7.3.2]{Illusie_Complexe_de_de_Rham_Witt_et_cohomologie_cristalline}, Illusie proves that the top Witt vector cohomology $H^2(S,W\cO_S)$ is $2$-torsion, and hence vanishes after twisting with $\QQ$.  It turns out (see \cite{Liedtke_Note_on_Non_Reduced_Picard_Schemes}) that in dimension two such examples arise only in characteristic two and three.  Unfortunately, we do not know if similar result holds in higher dimensions.
\end{example}

\begin{example}
\label{example:unpleasant_fibres}
We now present an example of a variety which is Gorenstein, normal, globally $F$-split and rational but fails to admit $\QQ_p$-rational quasi-resolution and hence luckily cannot be a general fibre of an MRCC fibration.  Let $E \subset \PP^2$ be a globally $F$-split (equiv. ordinary) elliptic curve.  We consider the pair $(\PP^2,E)$ which is log Calabi--Yau, that is, the divisor $K_{\PP^2} + E$ is linearly equivalent to zero.  Moreover, by \cite[Proposition 7.2]{Schwede_F_adjunction}, we know that the pair is also globally $F$-split.  Taking a blow-up $X$ of $\PP^2$ in ten points $P_i$ lying on $E$ we obtain a morphism of pairs $f \colon (X,\tilde E) \to (\PP^2,E)$, where $\tilde E$ is the strict transform of $E$.  Since the points lie on the curve the morphism is crepant, and therefore the pair $(X,\tilde E)$ is also globally $F$-split and log Calabi--Yau. Let $L_i$ be the exceptional curve over $P_i$.
The normal bundle of $\tilde E$ is anti-ample, and therefore there is hope for contracting it.   This can actually be performed using Keel's technique (see \cite{Keel_Basepoint_freeness_for_nef_and_big_line_bundle_in_positive_characteristics}) for the  big and nef divisor 
\[
\tilde E + \frac{1}{2}(L_1 + L_2),
\]
 assuming that 
\begin{equation}
\label{eq:unpleasant_fibres:torsion}
\left. \tilde E + \frac{1}{2}(L_0 + L_1) \right|_{\tE} \sim E|_E - \left(\sum_i P_i \right) + \frac{1}{2}(P_1 + P_2)
\end{equation}
is a torsion divisor on $E \cong \tE$. The points $P_i$ can always be chosen so  that the divisor of \autoref{eq:unpleasant_fibres:torsion} is torsion. So, let us assume that we made such a choice, and let $Y$ be the surface obtained by contracting $\tE$.  We claim that \emph{$Y$ satisfies our requirements -- it is Gorenstein, normal, globally $F$-split and rational}.  The canonical divisor $K_Y$ is the pushforward of $K_X + \tilde E$ and hence is trivial.  This implies that $Y$ is Gorenstein and normal.  The necessary $F$-splitting comes from a pushforward of a global $F$-splitting on $X$.  

The variety $Y$ is also interesting cohomologically.  More precisely, using proper base change and the fact that $E$ is ordinary we may prove that the natural resolution $X \to Y$ is not $\QQ_p$-rational. Then, the same holds for any other resolution too,  because $X \to Y$ is the minimal resolution.  Moreover, using Leray spectral sequence for W$\cO$-cohomology for the same resolution $f \colon X \to Y$ we obtain the following diagram describing degenerate $E_2$-page of the spectral sequence (we note that $f_*W\cO_{X,\QQ} = W\cO_{Y,\QQ}$):
\[
\xymatrix{
    H^0\left(Y,R^1f_*W\cO_{X,\bQ}\right) \expl{\isom}{\cite[Theorem 2.4]{Berthelot_Bloch_Esnault_On_Witt_vector_cohomology_for_singular_varieties}} H^1\left(E,W\cO_{E,\QQ}\right) \ar[rrd] & & \\ 
    H^0\left(Y,W\cO_{Y,\QQ}\right) & H^1\left(Y,W\cO_{Y,\QQ}\right) & H^2\left(Y,W\cO_{Y,\QQ}\right)
}
\]
which yields an isomorphism $H^2\left(Y,W\cO_{Y,\QQ}\right) \isom H^1\left(E,W\cO_{E,\QQ}\right)$ since $H^1\left(X,W\cO_{X,\QQ} \right) =H^2\left(X,W\cO_{X,\QQ} \right) = 0$ using rational chain connectedness.  As $E$ is ordinary, $0 \neq H^1\left(E,W\cO_{E,\QQ}\right) \cong H^1_{\et}(E,\QQ_p) \otimes_{\QQ_p} K$, and hence \autoref{lem:uniruled_witt_vector} does not hold true without appropriate singularity assumptions.
\end{example}